\documentclass[11pt,a4paper]{article}
\usepackage{amsmath,amssymb,amsthm}
\usepackage{amsfonts}

\usepackage{mathrsfs}
\usepackage{color}
\usepackage{ascmac}
\usepackage{dsfont}

\theoremstyle{definition}
 \newtheorem{dfn}{Definition}[section]

\theoremstyle{plain}
 \newtheorem{thm}[dfn]{Theorem}
 
 \newtheorem{lem}[dfn]{Lemma}

\newtheorem*{Rmk}{Remark}

\newcommand{\bn}{{\mathbf n}}

\newcommand{\bk}{{\mathbf k}}

\newcommand{\bu}{{\mathbf u}}

\newcommand{\bv}{{\mathbf v}}
\newcommand{\bw}{{\mathbf w}}

\newcommand{\bA}{{\mathbf A}}
\newcommand{\bB}{{\mathbf B}}

\newcommand{\bD}{{\mathbf D}}

\newcommand{\bI}{{\mathbf I}}

\newcommand{\bM}{{\mathbf M}}

\newcommand{\bV}{{\mathbf V}}
\newcommand{\bS}{{\mathbf S}}

\newcommand{\DV}{{\rm Div}\,}
\newcommand{\dv}{{\rm div}\,}

\newcommand{\BR}{{\mathbb R}}

\newcommand{\BN}{{\mathbb N}}

\newcommand{\CB}{{\mathcal B}}

\newcommand{\CD}{{\mathcal D}}
\newcommand{\CE}{{\mathcal E}}
\newcommand{\CF}{{\mathcal F}}
\newcommand{\CI}{{\mathcal I}}

\newcommand{\CM}{{\mathcal M}}

\newcommand{\CR}{{\mathcal R}}
\newcommand{\CS}{{\mathcal S}}

\newcommand{\CP}{{\mathcal P}}

\newcommand{\fa}{{\mathfrak a}}

\newcommand{\fh}{{\mathfrak h}}

\newcommand{\bg}{{\mathbf g}}
\newcommand{\bh}{{\mathbf h}}
\newcommand{\pd}{\partial}

\usepackage{xcolor}

\newcommand{\ep}{\varepsilon}

\newcommand{\wt}{\widetilde}

\newcommand{\ol}{\overline}

\newcommand{\loc}{{\rm loc\,}}

\numberwithin{equation}{section}

\usepackage{accents}

\newcommand{\Jt}{\langle t \rangle}
\newcommand{\Js}{\langle s \rangle}

\usepackage{scalerel,stackengine}
\stackMath
\newcommand\reallywidehat[1]{%
\savestack{\tmpbox}{\stretchto{%
  \scaleto{%
    \scalerel*[\widthof{\ensuremath{#1}}]{\kern-.6pt\bigwedge\kern-.6pt}%
    {\rule[-\textheight/2]{1ex}{\textheight}}
  }{\textheight}%
}{0.5ex}}%
\stackon[1pt]{#1}{\tmpbox}%
}

\newcommand{\vertiii}[1]{{\left\vert\kern-0.25ex\left\vert\kern-0.25ex\left\vert #1 
    \right\vert\kern-0.25ex\right\vert\kern-0.25ex\right\vert}}

\newcommand{\Ext}{{\rm Ext\,}}

\begin{document}

\title{\bf  Global wellposedness of the 3D compressible 
Navier-Stokes equations with free surface in the maximal regularity class}
\author{Yoshihiro Shibata
\thanks{Department of Mathematics,  Waseda University,
Ohkubo 3-4-1, Shinjuku-ku, Tokyo 169-8555, Japan
\endgraf
Department of Mechanical Engineering and Materials Science,
University of Pittsburgh, USA  \endgraf
e-mail address: yshibata325@gmail.com}
and Xin Zhang 
\thanks{School of Mathematical Sciences,
Tongji University, 
No.1239, Siping Road, Shanghai (200092), China. \endgraf
e-mail address: xinzhang2020@tongji.edu.cn
}}
\maketitle
\begin{abstract}
This paper concerns the global well posedness issue of the Navier-Stokes equations (CNS) describing barotropic compressible fluid flow 
with free surface occupied in the three dimensional exterior domain.  Combining
the maximal $L_p$-$L_q$ estimate and  the $L_p$-$L_q$ decay estimate of solutions 
to the linearized equations, we prove the unique existence of global in time 
solutions  in the time weighted  maximal $L_p$-$L_q$ regularity class for some $p>2$ and $q>3.$ Namely, the solution is bounded as $L_p$ in time and $L_q$ in space. Compared with the previous results of the free boundary value problem of (CNS) in unbounded domains,  we relax the regularity assumption on the initial states, which is the advantage by using the maximal $L_p$-$L_q$ regularity framework. On the other hand, the equilibrium state of the moving boundary of the exterior domain is not necessary the sphere.
To our knowledge, this paper is the first result on the long time solvability of the free boundary value problem of (CNS) in the exterior domain. 
\vskip1pc\noindent
subjclass[2020]: Primary: 35Q30; Secondary: 76N10 \vskip0.5pc\noindent
keywords: Compressible viscous barotropic fluid, free boundary problem, 
exterior domain, \\
\phantom{keywords } global well-posedness, $L_p$-$L_q$ maximal regularity, 
$L_p$-$L_q$ decay properties.
\end{abstract}


\section{Introduction}
\subsection{Model}
Let $\Omega_t$ be a time dependent exterior domain in the three dimensional
Euclidean space $\BR^3.$ To describe the motion of the viscous, barotropic,
 compressible fluid without surface tension, 
we denote the positive constant $\rho_e$ for the reference mass density,
$\Gamma_t$ for the moving boundary of $\Omega_t$ and $\bn_t$ for the unit outer normal to $\Gamma_t.$
The equations for the mass density $\rho+\rho_e>0,$ 
the velocity field $\bv=(v_1, v_2,  v_3)^\top,$ and $\Omega_t \subset \BR^3$ are formulated as follows:
\begin{equation}\label{eq:CNS}
	\left\{ \begin{aligned}
&\pd_t\rho+ \dv\big((\rho_e+\rho)\bv\big) =0
&&\quad&\text{in} \quad \bigcup_{0<t<T} \Omega_t \times\{t\}, \\
&(\rho_e+\rho)(\pd_t\bv +\bv\cdot\nabla\bv)
-\DV\big(\bS(\bv) - P(\rho_e+\rho)\bI\big) = 0 
&&\quad&\text{in} \quad \bigcup_{0<t<T} \Omega_t \times\{t\}, \\
&\big(\bS(\bv) - P(\rho_e+\rho)\bI\big)\bn_{\Gamma_t} = -P(\rho_e)\bn_{\Gamma_t},\,\,\,
V_{\Gamma_t} = \bv\cdot\bn_{\Gamma_t} 
&&\quad&\text{on} \quad \bigcup_{0<t<T} \Gamma_t \times\{t\}, \\
&(\rho, \bv,\Omega_t)|_{t=0} = (\rho_0, \bv_0,\Omega).
\end{aligned}
\right.
\end{equation}
Above, the Cauchy stress tensor 
$\bS(\bv) = \mu \bD(\bv) + (\nu-\mu)\dv \bv\bI 
\,\,\,\text{for the viscous coefficients}\,\, \mu,\nu>0,$
the doubled deformation tensor
 $\bD(\bv) = \nabla\bv + (\nabla\bv)^\top.$
Here, the $(i,j)$th entry of the matrix $\nabla \bv$ is $\pd_i v_j,$ 
$\bM^\top$ is the transposed  matrix of $\bM=[M_{ij}],$
the pressure law $P(\cdot)$ is a smooth function defined on $\BR_+,$ 
and $\bI$ is the $3\times 3$ identity matrix,
In addition, $\DV \bM$ denotes a $3$-vector of functions whose 
$i$-th component is $\sum_{j=1}^3 \pd_j M_{ij},$
$\dv \bv = \sum_{j=1}^3 \pd_jv_j,$ 
and $\bv\cdot\nabla = \sum_{j=1}^3 v_j\pd_j$ 
with $\pd_j = \pd/\pd x_j.$ 
In $\eqref{eq:CNS}_3,$ $P(\rho_e)$ coincides with the pressure of the atmosphere, and $V_{\Gamma_t}$ is the normal velocity of $\Gamma_t$
for $t>0.$
\medskip

In the mathematical theory of fluid mechanics, the barotropic (or isentropic) compressible Navier-Stokes system (CNS) is one fundamental model describing the dynamic behaviour of the viscous gases. There is a huge amount of literatures on the well posedness issue of (CNS). For the Cauchy problem of three dimensional viscous compressible flow, the global-in-time solution was first constructed by Matsumura and Nishida in \cite{MN1980} by assuming that the initial data are small enough in $H^3_2(\BR^3)$
\footnote{The notation on the functional space is given at the end of this section.}. 
 In addition, the authors \cite{MN1979} also proved that the solution
\footnote{$\rho$ denotes the perturbation of the mass density around some constant equilibrium state.}
 $(\rho,\bv)$ of (CNS) has the following decay property,
\begin{equation*}
\|(\rho,\bv)(\cdot,t)\|_{L_2(\BR^3)} \leq C (1+t)^{-3/4}
\end{equation*}
for some constant $C$ and $t>0,$ if the initial state $(\rho_0,\bv_0)$ is small measured by the norm $\|\cdot\|_{L_1(\BR^3) \cap H^3_2(\BR^3)}.$
Some improvement based on \cite{MN1979,MN1980} was obtained in the later works \cite{Kawa2002,LZ2011,WT2011}.
On the other hand, Hoff and Zumbrun in their seminal work \cite{HZ1995} studied the decay property of the solution in the $L_p$ norm, which exhibits the property of the multidimensional diffusion wave. For this aspect, the reader may also refer to \cite{KS1999,LW1998}. Furthermore, Danchin \cite{Dan2000} solved (CNS) with the $L_2$ Besov regularity. We also refer to \cite{CD2010,CMZ2010,Has2011} for the extension to the $L_p$ type Besov regularity and \cite{Okita2014,DanX2017} for the decay property. 
\smallbreak

In addition, the works \cite{ES2013,KK2005,KS2002,MN1983,ShiE2018} investigate
the non-slip boundary value problem of compressible flow in the smooth bounded domain, exterior domain or the half space. For the free boundary value problem in the smooth bounded domain without taking the surface tension into account, \cite{SV1983,Ta1981} solved the local wellposedness issue and \cite{Zaja1993,Shi2016c} established the global-in-time solution. When the moving bounded domain has the surface tension effect, \cite{SolTa1992,Zaja1994} proved the long time stability of some equilibrium state in the framework of the anisotropic Sobolev space. More recently,  \cite{EvBS2014,GS2014,Z2020} solved the free boundary value problem of (CNS) in short time interval but in the unbounded domain  within the maximal $L_p$-$L_q$ regularity framework.
\medskip

\subsection{Lagrangian coordinates and main result}
We assume that the reference (or initial) domain $\Omega$ is a $C^2$ exterior domain in $\BR^3.$ 
Since $\Omega_t$ is unknown, we shall transform $\Omega_t$ to the reference domain $\Omega$ by the Lagrangian transformation. 
  Let $\bu$ be a velocity field in the Lagrange 
coordinates $y = (y_1, y_2, y_3) \in \Omega$, and we define the Lagrangian mapping: 
\begin{equation}\label{def:LT}
x = X_{\bu}(y,t)= y + \int^t_0 \bu(y, s)\,ds
\end{equation}
for $0<t<T.$
In order to guarantee the invertibility of $X_{\bu}$ in \eqref{def:LT}, we assume that 
\begin{equation}\label{asmp:small_1}
\int^T_0\|\bu(\cdot, s)\|_{H^1_\infty(\Omega)}\,ds \leq \delta
\end{equation}
with small constant $\delta> 0.$ 
So we may write the inverse map of the transformation \eqref{def:LT} by $X_\bu^{-1}(x, t).$

Suppose that
\begin{gather*}
\rho(x, t) = \eta \big(X^{-1}_\bu(x. t), t\big), \quad
\bv(x, t) = \bu(X^{-1}_\bu(x, t), t),\\
\Omega_t = \{x =X_{\bu}(y,t): y \in \Omega\}, \quad
\Gamma_t = \{x =X_{\bu}(y,t) : y \in \Gamma=\pd \Omega\}.
\end{gather*}
The kinematic condition (non-slip condition):
$V_{\Gamma_t}=\bv\cdot\bn_t$ is automatically satisfied under the Lagrangian transformation
\eqref{def:LT}. 

Now, we derive the equations for $(\eta,\bu)$ 
from \eqref{eq:CNS}.
Introduce that
$$\frac{\pd x}{\pd y} =\nabla_y X_{\bu}= \bI + \int^t_0\nabla_y\bu(y, s)\,ds,$$
and $J_{\bu}=\det (\nabla_y X_{\bu}).$
Then by the assumption \eqref{asmp:small_1}, there exists the inverse of $\nabla_y X_{\bu},$ that is, 
$$\frac{\pd y}{\pd x} =\big(\nabla_y X_{\bu}\big)^{-1}
= \bI + \bV_0(\bk), \quad 
\bk =\int^t_0\nabla_y\bu(y, s)\,ds,
 $$
where $\bV_0(\bk)$ is a matrix-valued function given by the series
$$\bV_0(\bk)= \sum_{j=1}^{\infty} (-\bk)^j.$$
In particular, $\bV_0(0) = 0.$
By the chain  rule, we introduce the gradient, divergence and stress tensor operators with respect to the transformation \eqref{def:LT},
\begin{gather}
\nabla_{\bu} = \big(\bI + \bV_0(\bk) \big)\nabla_y,\quad
\dv_{\bu} \bu = \big(\bI+\bV_0(\bk)\big):\nabla_y \bu
= J_{\bu}^{-1}\dv_y \Big( J_{\bu} \big(\bI+ \bV_{0}(\bk) \big)^\top 
\bu\Big), \nonumber\\
\bD_{\bu} (\bu) =\big(\bI + \bV_0(\bk) \big)\nabla \bu 
+(\nabla \bu)^{\top}\big(\bI + \bV_0(\bk) \big)^{\top}
= \bD(\bu)+\bV_0(\bk) \nabla \bu
+\big(\bV_0(\bk) \nabla \bu\big)^{\top},
\label{sym:LL_1} \\
\bS_{\bu}(\bu)=\mu\bD_{\bu} (\bu)
 +(\nu-\mu) (\dv_{\bu} \bu) \bI,\quad 
\DV_{\bu} \bA= J_{\bu}^{-1} \DV_{y} \Big(J_{\bu}\bA \big(\bI + \bV_0(\bk) \big) \Big). \nonumber
\end{gather}
In addition, $\DV_{\bu} \bA$ can be also written via
\begin{equation}\label{eq:LL_DV_1}
\DV_{\bu} \bA= \DV_y \bA+\big(\bV_0(\bk) \nabla | \bA\big)
\end{equation}
with the $i$th component 
$(\bB \nabla | \bA)_{i}= \sum_{j,k=1}^3 B_{jk}\pd_k A_{ij}$ for $\bB=[B_{ij}]_{3\times 3},$ $i=1,2,3.$
Like the standard operator $\DV,$ $\DV_{\bu} \bA=0$ if $\bA$ is a constant matrix.
Then according to \eqref{sym:LL_1}, the pair $(\eta,\bu)$ fulfils the following system:
\begin{equation}\label{eq:LL_CNS_1}
	\left\{ \begin{aligned}
&\pd_t\eta + (\rho_e+\eta)\, \dv_{\bu} \bu=0
&&\quad&\text{in} &\quad \Omega \times (0,T), \\
&(\rho_e+\eta)\pd_t\bu-\DV_{\bu}\big(\bS_{\bu}(\bu) - P(\rho_e+\eta)\bI\big) = 0 
&&\quad&\text{in}& \quad \Omega \times (0,T), \\
&\big(\bS_{\bu}(\bu) - P(\rho_e+\eta)\bI\big)\bn_{\bu}
 = -P(\rho_e)\bn_{\bu}
&&\quad&\text{on}& \quad \Gamma \times (0,T), \\
&(\eta, \bu)|_{t=0} = (\rho_0, \bv_0)
&&\quad&\text{in}& \quad \Omega, \\
\end{aligned}
\right.
\end{equation}
where $\bn_{\bu}$ is defined by
$$\bn_{\bu} = \frac{\big(\bI + \bV_0(\bk)\big)\bn_{\Gamma}}{\big|\big(\bI + \bV_0(\bk)\big)\bn_{\Gamma} \big|} $$
with $\bn_{\Gamma}$ standing for the unit normal vector to $\Gamma=\pd \Omega.$ 
Our main result is the following global wellposedness for problem \eqref{eq:LL_CNS_1}.
\begin{thm}\label{thm:main_1} 
Let $\Omega$ be a $C^3$ exterior domain in $\BR^3.$ 
Let $b>0,$ $2<p<\infty$, $2<q_1<3 <q_2<\infty$ satisfying the conditions:
\begin{equation}\label{cdt:bq}
 \frac{1}{q_1}=\frac{1}{3}+\frac{1}{q_2},\quad  bp'>1,\quad
\frac{3}{2q_1}+\frac{1}{2} -\frac{1}{p}>b\geq \frac{3}{2q_1}\cdot
\end{equation}
Assume that $(\rho_0,\bv_0) \in\bigcap_{q=q_1, q_2}( H^1_q(\Omega) 
\times  B^{2-2/p}_{q,p}(\Omega)^{3}) \cap L_{q_1/2}(\Omega)^4$ satisfying the compatibility condition of order 1. 
Namely,
\begin{equation}\label{cdt:initial}
\big(\bS(\bv_0) - P(\rho_e+\rho_0)\bI\big)\bn_{\Gamma} = -P(\rho_e)\bn_{\Gamma}.
\end{equation}
Denote that 
\begin{equation*}
\CI_0 = \sum_{q\in \{q_1,q_2\}} 
\Big( \|\rho_0\|_{H^1_q(\Omega)} 
+ \|\bv_0\|_{B^{2-2/p}_{q,p}(\Omega)} \Big) 
+ \|(\rho_0, \bv_0)\|_{L_{q_1/2}(\Omega)}.
\end{equation*}
Then, there exists a small constant $\varepsilon > 0$ such that if 
$\CI_0 \leq \varepsilon,$ 
then problem \eqref{eq:LL_CNS_1} admits a unique solution $(\eta, \bu)$ with
\begin{align*}
\eta \in H^1_{p,{\rm loc}}\big(\BR_+; L_q(\Omega)\big), \quad
\bu  \in H^1_{p, {\rm loc}} \big(\BR_+; L_q(\Omega)^3\big) 
\cap L_{p, {\rm loc}}\big(\BR_+; H^2_q(\Omega)^3\big)
\end{align*}
for $q\in \{q_1,q_2\}.$
Moreover, there exists a constant $M$ independent of $T,$ $\varepsilon$ such that  
$\CE_T(\eta, \bu) \leq M\varepsilon $ for any $T > 0.$ 
Here, we have set 
\begin{equation}\label{def:CE_T}
\begin{aligned}
\CE_T(\eta, \bu)=
&\|\Jt^{3/(2q_1)}(\eta, \bu)\|_{L_\infty((0, T); L_{q_1}(\Omega))} 
+ \|\Jt^{b}(\eta, \bu)\|_{L_\infty((0, T); L_{q_2}(\Omega))} \\
&+\|\Jt^{b}( \nabla\eta,  \nabla\bu, \nabla^2 \bu)\|_{L_p((0,T);L_{q_1}(\Omega))}
+ \|\Jt^{b}  (\eta, \bu)\|_{L_p((0,T);H^1_{q_2}(\Omega))}  \\
 &+ \sum_{\ell=1,2}\Big(\|\Jt^{b} \pd_t\eta\|_{L_p((0,T);H^{1}_{q_\ell}(\Omega))}
+  \|\Jt^{b}  \pd_t\bu\|_{L_p((0,T); L_{q_\ell}(\Omega))} \Big)
\end{aligned}
\end{equation}
with $\Jt = \sqrt{1 + t^2}.$ 
\end{thm}

\begin{Rmk}Some comments on Theorem \ref{thm:main_1}:
\begin{enumerate}
\item In fact, it is not hard to see from \eqref{cdt:bq} in Theorem \ref{thm:main_1} that $q_2>6$ and $1/2<b<5/4.$  In particular, let $(q_1,q_2)=\big(2+2\sigma,6(1+\sigma)/(1-2\sigma)\big)$ for any small $\sigma>0,$ while the initial data have the $L_{1+\sigma}$ boundedness. Then we see 
the decay rate of $L_{2+2\sigma}$ norm of our solution from \eqref{def:CE_T}, namely,
\begin{equation}\label{eq:rmk}
\|(\eta,\bu)(\cdot,t)\|_{L_{2+2\sigma}(\Omega)} 
\leq C_0 (1+t)^{-\frac{3}{4(1+\sigma)}}
\end{equation}
for some constant $C_0$ depending on $\CI_0.$
Moreover, \eqref{eq:rmk} is optimal in sense of heat flows.

\item Although $\bv_0$ does not necessarily belong to $B^{2-2/p}_{q_1/2,p}(\Omega),$ the maximal $L_{p}$-$L_{q_1/2}$ regularity property is used in later proof (see Theorem \ref{thm:nonlinear}). So we cannot extend the result to $q_1=1$ by the method of this paper, which distincts our result with the classical works \cite{MN1979,MN1980} in the $L_2$ energy approach.

\item Notice that $(\rho_e,0,\Omega)$ can be regarded as the equilibrium state of the problem \eqref{eq:CNS}. Furthermore, only $C^3$ smoothness constrain is imposed on $\Omega,$ but not the shape of the initial boundary $\Gamma.$ This is the main difference between the surface tension problem and the problem of this paper.

\item For technical reason, we only construct the global solution of \eqref{eq:LL_CNS_1} with the constrain $p>2$ in Lagrangian coordinates. However, we also expect the long time result for the endpoint case $p=2$ in the future.
\end{enumerate}

\end{Rmk}

\subsection{Reformulation of \eqref{eq:LL_CNS_1} and main idea}
Now let us give some rough idea of the proof of Theorem \ref{thm:main_1}.
To convert \eqref{eq:LL_CNS_1} into some linearized form, we set some notations.
It is clear that the boundary condition in \eqref{eq:LL_CNS_1} is equivalent to
\begin{equation}\label{eq:bc_1}
\Big(\bS_{\bu}(\bu) - \big(P(\rho_e+\eta)-P(\rho_e)\big)\bI\Big)\big(\bI + \bV_0(\bk)\big)\bn_{\Gamma} = 0.
\end{equation}
On the other hand, Taylor's theorem yields that 
\begin{equation*}
P(\rho_e+\eta) - P(\rho_e) = P'(\rho_e)\eta +Q(\eta) 
\,\,\,\text{with}\,\,\, 
Q(\eta)=\eta^2 \int_0^1 P''(\rho_e +\theta\eta) (1-\theta) \,d\theta.
\end{equation*}
From \eqref{sym:LL_1}, we introduce that 
\begin{equation}\label{sym:LL_2}
\begin{aligned}
\CD_{\bD}(\bu)&= \bD_{\bu} (\bu) -\bD(\bu)
= \bV_0(\bk) \nabla \bu
+\big(\bV_0(\bk) \nabla \bu\big)^{\top},\\
\CS_{\bD}(\bu) &= \bS_{\bu} (\bu) -\bS(\bu)
=\mu \CD_{\bD}(\bu) 
+ (\nu-\mu) \big( \bV_0(\bk): \nabla \bu \big) \bI.
\end{aligned}
\end{equation}
Suppose that $\Omega^c \subset B_R$ for some $R>0.$ 
Let $\kappa \in C_0^{\infty}(\BR^N)$ be a cut-off function such that 
$\kappa (x)=1$ for $|x|\leq 2R$ and $\kappa(x)=0$ for $|x|\geq 3R.$ 
Then we rewrite the system \eqref{eq:LL_CNS_1} in the following equivalent form:
\begin{equation}\label{eq:LL_CNS_2}
	\left\{ \begin{aligned}
&\pd_t\eta+ \gamma_1 \, \dv \bu=f(\eta,\bu)
&&\quad&\text{in} &\quad \Omega \times (0,T), \\
&\gamma_1 \pd_t\bu 
-\DV\big(\bS(\bu) - \gamma_2 \eta\bI\big) = \bg(\eta,\bu) 
&&\quad&\text{in}& \quad \Omega \times (0,T), \\
&\big(\bS(\bu) - \gamma_2 \eta \bI\big)\bn_{\Gamma} = \kappa\bh(\eta,\bu)
&&\quad&\text{on}& \quad \Gamma \times (0,T), \\
&(\eta, \bu)|_{t=0} = (\rho_0, \bv_0)
&&\quad&\text{in}& \quad \Omega, \\
\end{aligned}
\right.
\end{equation}
where $\gamma_1= \rho_e,$ $\gamma_2=P'(\rho_e)$ and 
the nonlinear terms on the right-hand side are given by
\begin{equation}\label{def:fgh}
\begin{aligned}
f(\eta,\bu)=& -\eta \,\dv \bu -(\gamma_1+\eta) \bV_0(\bk):\nabla \bu,\\
\bg(\eta,\bu)=& -\eta \pd_t \bu 
+\DV\big(\CS_{\bD}(\bu) - Q(\eta)\bI\big)
+\big( \bV_0(\bk) \nabla \mid \bS_{\bu}(\bu) - \gamma_2 \eta \bI - Q(\eta)\bI\big),\\
\bh(\eta,\bu)=&-\big(\CS_{\bD}(\bu)-Q(\eta)\bI\big) \bn_{\Gamma}
-\big(\bS_{\bu}(\bu) - \gamma_2 \eta \bI -Q(\eta) \bI\big) 
\bV_0(\bk) \bn_{\Gamma}.
\end{aligned}
\end{equation}

Now we decompose the solution $(\eta,\bu)$ of \eqref{eq:LL_CNS_2} by 
$(\eta,\bu)=(\rho^1,\bv^1)+(\rho^2, \bv^2) + (\theta,\bw)$
where $(\rho^i,\bv^i),$ $i=1,2,$ 
satisfy the shifted systems for some $\lambda_0>0$: 
\begin{equation} \label{eq:sft1}
	\left\{ \begin{aligned}
&\pd_t\rho^1 + \lambda_0 \rho^1+ \gamma_1 \, \dv\bv^1= f(\eta,\bu)
&&\quad&\text{in} &\quad \Omega \times (0,T), \\
&\gamma_1 \pd_t\bv^1 +\lambda_0 \bv^1
-\DV\big(\bS( \bv^1 ) - \gamma_2 \rho^1 \bI\big)
= \bg(\eta,\bu)
&&\quad&\text{in}& \quad \Omega  \times (0,T), \\
&\big(\bS(\bv^1) - \gamma_2 \rho^1 \bI\big)\bn_{\Gamma} =\kappa \bh(\eta,\bu)
&&\quad&\text{on}& \quad \Gamma \times (0,T); 
&&\quad&\text{in}& \quad \Omega; 
\end{aligned}
\right.
\end{equation}
\begin{equation} \label{eq:sft2}
	\left\{ \begin{aligned}
&\pd_t\rho^2 + \lambda_0 \rho^2+ \gamma_2 \, \dv\bv^2= \lambda_0\rho^1
&&\quad&\text{in} &\quad \Omega \times (0,T), \\
&\gamma_1 \pd_t\bv^2 +\lambda_0 \bv^2
-\DV\big(\bS( \bv^2 ) - \gamma_2 \rho^2 \bI\big)= \lambda_0\bv^1
&&\quad&\text{in}& \quad \Omega  \times (0,T), \\
&\big(\bS(\bv^2) - \gamma_2 \rho^2 \bI\big)\bn_{\Gamma} =0
&&\quad&\text{on}& \quad \Gamma \times (0,T), \\
&(\rho^2, \bv^2)|_{t=0} = (0, 0)
&&\quad&\text{in}& \quad \Omega.
\end{aligned}
\right.
\end{equation}
In addition, $(\theta, \bw)$ satisfies the compensate system: 
\begin{equation} \label{eq:homo}
	\left\{ \begin{aligned}
&\pd_t\theta + \gamma_1 \, \dv\bw=  \lambda_0 \rho^2 
&&\quad&\text{in} &\quad \Omega \times (0,T), \\
&\gamma_1 \pd_t\bw -\DV\big(\bS(\bw) - \gamma_2 \theta \bI\big)
=\lambda_0 \bv^2 
&&\quad&\text{in}& \quad \Omega  \times (0,T), \\
&\big(\bS(\bw) - \gamma_2 \theta \bI\big)\bn_{\Gamma} =0
&&\quad&\text{on}& \quad \Gamma \times (0,T), \\
&(\theta, \bw )|_{t=0} =(\rho_0, \bv_0) - (\rho^1, \bv^1)|_{t=0}
&&\quad&\text{in}& \quad \Omega.
\end{aligned}
\right.
\end{equation}

In the system \eqref{eq:sft1},  we omit the initial conditions. 
Because, the role of $(\rho^1,\bv^1)$ is to handle the nonlinear terms as in Section \ref{sec:shift}. To apply the semigroup theory established in \cite{ShiZ2020a}, the construction of $(\rho^2,\bv^2)$ is necessary, which can be regarded as one difficulty from the quasilinear type conditions \eqref{eq:bc_1} in our problem.
At last, we treat \eqref{eq:homo} by the decay estimates of the semigroup theory (see Sections \ref{sec:mi} and \ref{sec:decay}).

\subsection{Notion}\label{subsec:notion}
Let us end up this section by introducing some symbols for the domains and the functional spaces, which will be used throughout this paper.
For any domain $G \subset \BR^3,$ $1\leq p \leq \infty$ and $k \in \BN_0=\BN \cup \{0\},$ 
$L_p(G)$ and $H^k_p(G)$ denote the standard Lebesgue space and Sobolev space respectively. 
In addition, the Besov space $B^{s}_{q,p}(G)$ 
for some $k-1<s\leq k$ and $(p,q) \in (1,\infty)^2$
is defined by the real interpolation functor
\begin{equation*}
B_{q,p}^{s}(G)= \big(L_q(G),H^{k}_q(G)\big)_{s\slash k,p}.
\end{equation*} 
For any Banach space $X$ and any interval $I$ in $\BR,$ $H^k_p(I;X)$ (or $H^k_{p,\loc}(I;X)$)  stands for
the total of the $X$-valued mappings in the $H^k_p$ (or $H^k_{p,\loc}$) class.
Sometimes, we write $H^k_p(a,b;X)$ for simplicity  whenever $I=(a,b)$ for any 
$0\leq a<b<\infty.$
In addition, $A\lesssim B$ stands for $A\leq C B$ whenever $C$ is a harmless constant.

%


\section{Linear theory on the shifted model}
In this section, we introduce the linear theory on the time-weight estimates of  some shifted model problem.
\begin{thm}\label{thm:linear}
Let $\Omega$ be a $C^2$ exterior domain in $\BR^3,$ and let $1<p,q<\infty.$  
Assume that $f,$ $\bg$ and $\bh$ satisfy that 
\begin{gather*}
 f (\cdot,t) \in L_{p}\big(0, T;H^1_q(\Omega)\big), \quad
\bg (\cdot,t) \in L_{p}\big(0, T ;L_q(\Omega)^3\big), \\
\Jt^{b} \bh (\cdot,t) \in L_{p}\big(\BR;H^1_q(\Omega)^3\big) \cap 
H^{1/2}_{p}\big(\BR;L_q(\Omega)^3\big)
\end{gather*}  
for some $b\geq 0.$ 
Let $f_0$ and $\bg_0$ be zero extensions of $f$ and $\bg$ outside of $(0, T)$, that is
$h_0(\cdot, t) = h(\cdot, t)$ for $t \in (0, T)$ and $h_0(\cdot, t) = 0$ for 
$t \not\in (0, T)$ with $h \in \{f, \bg\}$. 
Then there exists a constant $\lambda_0>0$ such that the following system 
\begin{equation}\label{shift.eq.1}
	\left\{ \begin{aligned}
&\pd_t\rho + \lambda_0 \rho+ \gamma_1 \, \dv\bv= f_0
&&\quad&\text{in} &\quad \Omega \times \BR, \\
&\gamma_1 \pd_t\bv + \lambda_0 \bv
-\DV\big(\bS( \bv ) - \gamma_2 \rho \bI\big)= \bg_0
&&\quad&\text{in}& \quad \Omega  \times \BR, \\
&\big(\bS(\bv) - \gamma_2 \rho \bI\big)\bn_{\Gamma} = \bh
&&\quad&\text{on}& \quad \Gamma \times \BR
\end{aligned}
\right.
\end{equation}
admits a unique solution $(\rho, \bv)$ with
$$\rho \in H^1_{p}\big(\BR; H^1_q(\Omega)\big), 
\quad \bv \in H^1_{p}\big(\BR; L_q(\Omega)^3\big) 
\cap L_{p}\big(\BR; H^2_q(\Omega)^3\big),$$
possessing the estimate:
\begin{equation}\label{es:Lame_1}
\begin{aligned}
& \|\Jt^{b}(\rho, \pd_t \rho)\|_{L_p((0, T); H^1_q(\Omega))}
+ \|\Jt^{b}\pd_t\bv \|_{L_p((0, T); L_q(\Omega))} 
+\|\Jt^{b} \bv\|_{L_p((0, T); H^2_q(\Omega))} \\
\leq & C \Big(
\|\Jt^{b} f\|_{L_p((0, T); H^1_q(\Omega))}
+ \|\Jt^{b}\bg\|_{L_p((0, T); L_q(\Omega))}
+ \|\Jt^{b}\bh\|_{L_p(\BR; H^1_q(\Omega))} 
+ \|\Jt^{b}\bh\|_{H^{1/2}_p(\BR; L_q(\Omega))}\Big)
\end{aligned}
\end{equation}
for some constant $C > 0.$
\smallbreak

Moreover, if $\bh=0$, then $(\rho, \bv) = (0, 0)$ for $t \leq 0$, in particular
$(\rho, \bv)|_{t=0} = (0, 0)$.
\end{thm}

\begin{proof}
One can establish \eqref{es:Lame_1} for the case $b=0$ by using the
$\CR$ bounded solution operators of the corresponding generalized
resolent problem  in \cite[Theorem 2.4]{EvBS2014} 
(see also \cite[Section 5]{ShiZ2020a}). Namely, we have for some $\lambda_0>0$
\begin{equation}\label{es:Lame_2}
\begin{aligned}
& \|(\rho, \pd_t \rho)\|_{L_p(\BR; H^1_q(\Omega))}
+ \|\pd_t\bv \|_{L_p(\BR; L_q(\Omega))} 
+\|\bv\|_{L_p(\BR, H^2_q(\Omega))} \\
\leq & C \Big(
\|f_0\|_{L_p(\BR; H^1_q(\Omega))}
+ \|\bg_0\|_{L_p(\BR; L_q(\Omega))} 
+ \|\bh\|_{L_p(\BR; H^1_q(\Omega))}
+ \|\bh\|_{H^{1/2}_p(\BR; L_q(\Omega))}\Big).
\end{aligned}
\end{equation}
Then we can deduce the bound \eqref{es:Lame_1} for general $b>0$ 
by considering the following system
\begin{equation} \label{eq:Linear_tb}
	\left\{ \begin{aligned}
&\pd_t (\Jt^b \rho) + \lambda_0 \Jt^b \rho + \gamma_1 \, \dv (\Jt^b \bv)= \wt f
&&\quad&\text{in} &\quad \Omega \times \BR, \\
&\gamma_1 \pd_t(\Jt^b\bv) + \lambda_0 \Jt^b\bv
-\DV\big(\bS( \Jt^b\bv ) - \gamma_2 \Jt^b \rho \bI\big)= \wt \bg
&&\quad&\text{in}& \quad \Omega  \times \BR, \\
&\big(\bS(\Jt^b\bv) - \gamma_2 \Jt^b \rho \bI\big)\bn_{\Gamma} = \wt \bh
&&\quad&\text{on}& \quad \Gamma \times \BR,
\end{aligned}
\right.
\end{equation}
where we have set $\wt{f}=\Jt^b f + b \Jt^{b-2} t \rho,$
$\wt{\bg}= \Jt^b \bg +\gamma_1 b \Jt^{b-2}t \bv$ and
$\wt{\bh}=\Jt^b \bh.$ 
Note the fact that 
\begin{equation}\label{es:Jt_1}
\Jt^{b-1} \leq 
\begin{cases}
1& \text{for}\,\,\, 0<b\leq 1,\\
C_{\delta} + \delta \Jt^{b}
& \text{for}\,\,\, b>1\,\,\, \text{and any}\,\,\,  0<\delta<1.
\end{cases}
\end{equation}
Thus we see from \eqref{es:Lame_2} and \eqref{es:Jt_1} that 
\begin{align*}
&\|\wt f_0\|_{L_p(\BR; H^1_q(\Omega))}
+ \|\wt \bg\|_{L_p(\BR; L_q(\Omega))}
+ \|\wt \bh\|_{H^{1/2}_p(\BR; L_q(\Omega))}
+ \|\wt \bh\|_{L_p(\BR; H^1_q(\Omega))} \\
\leq& C_{\delta}\Big(
 \|\Jt^{b}f\|_{L_p(\BR; H^1_q(\Omega))}
+ \|\Jt^{b}\bg\|_{L_p(\BR; L_q(\Omega))}
+\|\Jt^{b} \bh\|_{L_p(\BR; H^1_q(\Omega))} 
 + \|\Jt^{b}\bh\|_{H^{1/2}_p(\BR; L_q(\Omega))} \Big)\\
&+\delta \big( \|\Jt^{b}\rho\|_{L_p(\BR; H^1_q(\Omega))} 
+\|\Jt^{b}\bv\|_{L_p(\BR; L_q(\Omega))} \big).
\end{align*}
Then applying \eqref{es:Lame_2} to \eqref{eq:Linear_tb} and choosing $\delta$ small yield the desired bound of $(\rho,\bv).$
\medskip

Moreover, if $\bh=0$, then for any $\gamma >\lambda_0$ we see from \cite[Theorem 2.4]{EvBS2014} that 
\begin{equation}
\begin{aligned}
\gamma \|(\rho,\bv)\|_{L_p((-\infty,0]; L_q(\Omega))}
\leq & \gamma \|e^{-\gamma t}(\rho,\bv)\|_{L_p(\BR; L_q(\Omega))}  \\
\leq & C \Big(
\|e^{-\gamma t} f_0\|_{L_p(\BR; H^1_q(\Omega))}
+ \|e^{-\gamma t} \bg_0\|_{L_p(\BR; L_q(\Omega))} \Big)\\
\leq &  C \Big(
\|f\|_{L_p(0,T; H^1_q(\Omega))}
+ \|\bg\|_{L_p(0,T; L_q(\Omega))} \Big).
\end{aligned}
\end{equation}
Letting $\gamma \to\infty$, we have $(\rho, \bv) = (0, 0)$ for $t \leq 0.$
\smallbreak

The uniqueness of the solution follows from the uniqueness of the 
corresponding resolvent problem after applying the Laplace transform
to equations \eqref{shift.eq.1}. 
\end{proof}

\if
In fact, Theorem \ref{thm:linear} can be slightly refined to study the local-in-time model as follows
\begin{equation}\label{eq:linear_shift}
	\left\{ \begin{aligned}
&\pd_t\rho + \lambda_0 \rho+ \gamma_1 \, \dv\bv= f
&&\quad&\text{in} &\quad \Omega \times (0,T), \\
&\gamma_1 \pd_t\bv + \lambda_0 \bv
-\DV\big(\bS( \bv ) - \gamma_2 \rho \bI\big)= \bg
&&\quad&\text{in}& \quad \Omega  \times (0,T), \\
&\big(\bS(\bv) - \gamma_2 \rho \bI\big)\bn_{\Gamma} = \bh
&&\quad&\text{on}& \quad \Gamma \times (0,T), \\
&(\rho, \bv )|_{t=0} =(\rho_0,\bv_0)
&&\quad&\text{in}& \quad \Omega.
\end{aligned}
\right.
\end{equation}
To this end, we introduce some time-weighted extension operator which is also important for later study.
For any $b\geq 0,$ let us set
\begin{equation*}
E_{(T)}^b \fh(\cdot, t) = 
\begin{cases}
		\fh(\cdot, t) &\quad\mbox{if}\quad 0<t<T, \\
		\frac{\langle 2T-t \rangle^b}{\Jt^b}\fh(\cdot, 2T-t) &\quad\mbox{if}\quad  T<t<2T, \\
		0 &\quad\mbox{otherwise}.
	\end{cases}
\end{equation*}
Here we give some result concerning the extension operator $E_{(T)}^b.$
\begin{lem}
Let $X$ be a Banach space, $b\geq 0$ and $1\leq r\leq \infty.$
Assume that $\fh$ is some $X$-valued function. Then we have 
\begin{equation}\label{es:Fh_1}
\|\Jt^b E_{(T)}^b\fh(\cdot, t) \|_{L_{r} (\BR; X)} 
\leq 2^{1/r} \|\Jt^b \fh(\cdot, t) \|_{L_r (0,T; X)}.
\end{equation}
In addition, if $\fh(\cdot,0)=0,$ then there holds 
\begin{equation}\label{es:Fh_2}
\|\Jt^b \pd_t E_{(T)}^b\fh(\cdot, t) \|_{L_{r} (\BR; X)} 
\leq  C\big(\|\Jt^b \fh(\cdot, t) \|_{L_r (0,T; X)} 
+\|\Jt^b \pd_t \fh(\cdot, t) \|_{L_r (0,T; X)} \big)
\end{equation}
for some constant $C.$
\end{lem}
\begin{proof}
By the definition of $E_{(T)}^b,$ we immediately have \eqref{es:Fh_1} and 
\begin{equation}\label{es:Fh_0}
\|\pd_t E_{(T)}^0\fh(\cdot, t) \|_{L_{r} (\BR; X)} 
\leq 2^{1/r} \|\pd_t \fh(\cdot, t) \|_{L_r (0,T; X)}.
\end{equation}
Let us consider \eqref{es:Fh_2} for $b>0.$ Notice that $\langle 2T-t \rangle \leq  \Jt$ for  $T \leq t \leq 2T.$ Then \eqref{es:Jt_1} implies that
\begin{equation}\label{es:Jt_2}
\frac{d}{dt} \Big( \frac{\langle 2T-t\rangle^b}{\Jt^b} \Big) 
\leq  2b \frac{\langle 2T-t\rangle^{b-1}}{\Jt^b}  \leq 
\begin{cases}
2b \Jt^{-b}  & (0<b\leq 1),\\
C_{b,\delta} \Jt^{-b} + \delta \langle 2T-t\rangle^b \Jt^{-b}
&  (b>1),
\end{cases}
\end{equation}
for any  $T \leq t \leq 2T$ and $0<\delta<1.$ 
Thus, it is not hard to show \eqref{es:Fh_2} from \eqref{es:Jt_2}.

\end{proof}

By using the extension operator $E_{(T)}^b,$ we have the following result.
\begin{thm}\label{thm:linear_2}
Let $\Omega$ be a $C^2$ exterior domain in $\BR^3,$ and let $1<p,q<\infty.$  
Assume that $f,$ $\bg$ and $\bh$ satisfy that 
\begin{gather*}
\Jt^{b} f (\cdot,t) \in L_{p}\big(0,T;H^1_q(\Omega)\big), \quad
\Jt^{b} \bg (\cdot,t) \in L_{p}\big(0,T;L_q(\Omega)^3\big), \\
\Jt^{b} \bh (\cdot,t) \in H^{1/2}_{p}\big(\BR;L_q(\Omega)^3\big)
\cap  L_{p}\big(\BR;H^1_q(\Omega)^3\big)
\end{gather*}  
for some $b\geq 0.$ Let $(\rho_0, \bv_0)$ belong to 
$H^1_q(\Omega)\times B^{2-2/p}_{q,p}(\Omega)^{3}$ fulfilling
the compatibility conditions:
\begin{equation*}
\big(\bS(\bv_0) - \gamma_2 \rho_0 \bI\big)\bn_{\Gamma} = \bh|_{t=0+}.
\end{equation*}
Then there exists a constant $\lambda_0>0$ such that 
the system \eqref{eq:linear_shift} admits a unique solutions $(\rho, \bv)$ with
$$\rho \in H^1_{p}\big(0,T; H^1_q(\Omega)\big), 
\quad \bv \in H^1_{p}\big(0,T; L_q(\Omega)^3\big) 
\cap L_{p}\big(0,T; H^2_q(\Omega)^3\big),$$
possessing the estimate:
\begin{equation}\label{es:Lame_T}
\begin{aligned}
& \|\Jt^{b}(\rho, \pd_t \rho)\|_{L_p(0,T; H^1_q(\Omega))}
+ \|\Jt^{b}\pd_t\bv \|_{L_p(0,T; L_q(\Omega))} 
+\|\Jt^{b} \bv\|_{L_p(0,T; H^2_q(\Omega))} \\
\leq & C\Big(\|\rho_0\|_{H^1_q(\Omega)} + \|\bv_0\|_{B^{2-2/p}_{q,p}(\Omega)} 
+ \|\Jt^{b} f\|_{L_p(0,T; H^1_q(\Omega))}
+ \|\Jt^{b}\bg\|_{L_p(0,T; L_q(\Omega))}\\
&\quad \quad  
+ \|\Jt^{b}\bh\|_{H^{1/2}_p(\BR; L_q(\Omega))}
+ \|\Jt^{b}\bh\|_{L_p(\BR; H^1_q(\Omega))} \Big)
\end{aligned}
\end{equation}
for some constant $C > 0.$
\end{thm}
\begin{proof}
By applying \eqref{es:Lame_1} and \eqref{es:Fh_1} to the following system 
\begin{equation*}
	\left\{ \begin{aligned}
&\pd_t\wt \rho + \lambda_0\wt \rho+ \gamma_1 \, \dv \wt \bv=E_{(T)}^b f
&&\quad&\text{in} &\quad \Omega \times \BR_+, \\
&\gamma_1 \pd_t \wt \bv + \lambda_0 \wt \bv
-\DV\big(\bS( \wt \bv ) - \gamma_2\wt \rho\, \bI\big)=E_{(T)}^b \bg
&&\quad&\text{in}& \quad \Omega  \times \BR_+, \\
&\big(\bS(\wt \bv) - \gamma_2 \wt \rho \,\bI\big)\bn_{\Gamma} = \bh
&&\quad&\text{on}& \quad \Gamma \times \BR_+, \\
&(\wt \rho, \wt\bv )|_{t=0} =(\rho_0,\bv_0)
&&\quad&\text{in}& \quad \Omega,
\end{aligned}
\right.
\end{equation*}
it is not hard to gain \eqref{es:Lame_T}. 
This completes the proof.
\end{proof}

\fi

\section{The solution of the shifted system}\label{sec:shift}
In this whole section,  we set for convenience that 
\begin{equation*}
\CM^{b}_{p,q}(T;\rho,\bv)=
 \|\Jt^{b}(\rho, \pd_t \rho)\|_{L_p(0,T; H^1_q(\Omega))}
+ \|\Jt^{b}\pd_t\bv \|_{L_p(0,T; L_q(\Omega))} 
+\|\Jt^{b} \bv\|_{L_p(0,T; H^2_q(\Omega))},
\end{equation*}
for any $b\geq 0,$ $1< p,q <\infty$ and $0<T\leq \infty.$ For the solution $(\eta,\bu)$ of \eqref{eq:LL_CNS_1}, recall that
\begin{equation} \label{def:X_T}
\begin{aligned}
 \CE_T (\eta,\bu) =
  &\|\Jt^{3/(2q_1)}(\eta, \bu)\|_{L_\infty(0, T; L_{q_1}(\Omega))}
 +\|\Jt^{b} \nabla \eta \|_{L_p(0,T;L_{q_1}(\Omega))}
 + \|\Jt^{b} \pd_t \eta\|_{L_p(0,T;H^1_{q_1}(\Omega))} \\
 &+\|\Jt^{b}  \pd_t \bu\|_{L_p(0,T;L_{q_1}(\Omega))}  
+ \|\Jt^{b}  \nabla \bu\|_{L_p(0,T;H^1_{q_1}(\Omega))} \\
&+ \|\Jt^b(\eta, \bu)\|_{L_\infty(0, T; L_{q_2}(\Omega))} 
+ \CM^{b}_{p,q_{2}}(T;\eta,\bu)\\
\end{aligned}
\end{equation}
for $b>0,$ $1<p,q_1,q_2<\infty$ and $0<T<\infty.$
We consider shifted equations:
\begin{equation} \label{eq:shift}
	\left\{ \begin{aligned}
&\pd_t\rho^1 + \lambda_0 \rho^1+ \gamma_1 \, \dv\bv^1= f(\eta,\bu)_0
&&\quad&\text{in} &\quad \Omega \times \BR, \\
&\gamma_1 \pd_t\bv^1 +\lambda_0 \bv^1
-\DV\big(\bS( \bv^1 ) - \gamma_2 \rho^1 \bI\big)= \bg(\eta,\bu)_0
&&\quad&\text{in}& \quad \Omega  \times \BR, \\
&\big(\bS(\bv^1) - \gamma_2 \rho^1 \bI\big)\bn_{\Gamma} 
=\Jt^{-b}\kappa \Ext[\bh(\eta,\bu)]
&&\quad&\text{on}& \quad \Gamma \times \BR.
\end{aligned}
\right.
\end{equation}
Here, $\Ext[\bh(\eta, \bu)]$ is an extension of $\bh(\eta,\bu)$ which will be given
later. 
In what follows, we study the system \eqref{eq:shift} and 
prove: 
\begin{thm} \label{thm:nonlinear}
Let $b>0,$ $1<p<\infty,$ $2<q_1\leq q_2$ and $3<q_2<\infty$ with $bp'>1.$
In addition, assume that $(\eta,\bu)$ satisfies the conditions:
\begin{equation}\label{cdt:etau} 
\CE_T=\CE_T (\eta,\bu)<\infty,\quad
\|\eta\|_{L_\infty(0,T;L_{\infty}(\Omega))}  \leq \rho_e/2, \quad 
\int^T_0\|\nabla\bu(\cdot, s)\|_{L_{\infty}(\Omega)}\,ds\leq 1/2.
\end{equation}
Let $\rho_0 \in H^1_{q_1}(\Omega) \cap H^1_{q_2}(\Omega)$ in \eqref{eq:LL_CNS_2}. 
Then \eqref{eq:shift} admits a unique solution $(\rho^1,\bv^1)$ with 
\begin{equation*}
\rho^1 \in H^1_p(\BR;H^1_q(\Omega))
\,\,\,\hbox{and}\,\,\, 
\bv^1 \in H^1_p(\BR;L_q(\Omega)^3) 
\cap L_p(\BR;H^2_q(\Omega)^3)
\end{equation*}
for any $q\in \{q_1/2,\, q_1, q_2\}.$
Moreover, there exists a constant $C$ such that the following estimate hold:
\begin{equation}\label{es:nonlinear}
\sum_{q\in \{q_1/2,\, q_1, q_2\}}  \CM^{b}_{p,q}(T;\rho^1,\bv^1)
\leq  C (\CI_0' +\CE_T^2 +\CE_T^3 ).
\end{equation}
Here, $\CI_0' = \sum_{\ell=1, 2} \|\rho_0\|_{H^1_{q_\ell}(\Omega)}$.
\end{thm}

The rest of this section is dedicated to studying the quantities on the right-hand side of \eqref{eq:shift}. Without loss of generality, we assume that $\CI_0' \leq 1$ hereafter. 

\subsection{Bound of $f(\eta,\bu)$}
We review some technical lemma proved in \cite[Appendix A]{SSZ2020}, which is about the estimates of $\bV_0(\bk).$
\begin{lem}[\cite{SSZ2020}] \label{lem:es_V}
Assume that $\bu$ is some smooth enough vector field satisfying 
\begin{equation*}
\|\nabla \bu\|_{L_1(0,T;L_{\infty}(\Omega))} \leq \sigma <1
\end{equation*}
for some constant $\sigma.$ Then there exists a positive constant $C_{\kappa}$ such that 
\begin{equation*}
\begin{aligned}
\|\bV_0(\bk)\|_{L_{\infty}(0,T;L_{\infty}(\Omega))} 
& \leq C_{\sigma} \|\nabla \bu\|_{L_1(0,T;L_{\infty}(\Omega))}, \\
\|\nabla \bV_0(\bk)\|_{L_{\infty}(0,T;L_{q}(\Omega))} 
& \leq C_{\sigma} \|\nabla^2 \bu\|_{L_1(0,T;L_{q}(\Omega))},\\
\|\pd_t \bV_0(\bk)\|_{L_{p}(0,T;H^1_{q}(\Omega))} 
& \leq C_{\sigma} \|\nabla \bu\|_{L_p(0,T;H^1_{q}(\Omega))}
\end{aligned}
\end{equation*}
for $1<q<\infty.$
\end{lem}
First, we notice from Lemma \ref{lem:es_V},  Sobolev inequalities and the condition $bp'>1$ that 
\begin{equation}\label{es:V_1}
\begin{aligned}
\|\bV_0(\bk)\|_{L_{\infty}(0,T;L_{\infty}(\Omega))} 
&\lesssim \|\nabla \bu\|_{L_1(0,T;H^{1}_{q_2}(\Omega))}  
\lesssim \|\Jt^b \nabla \bu\|_{L_p(0,T;H^{1}_{q_2}(\Omega))},\\
\|\nabla \bV_0(\bk)\|_{L_{\infty}(0,T;L_{q_\ell}(\Omega))}
&\lesssim \|\nabla^2 \bu\|_{L_1(0,T;L_{q_\ell}(\Omega))}  
\lesssim \|\Jt^b \nabla^2 \bu\|_{L_p(0,T;L_{q_\ell}(\Omega))}
\end{aligned}
\end{equation}
for $\ell=1,2.$
On the other hand, assuming the condition 
$\|\nabla \bu\|_{L_1(0,T;L_{\infty}(\Omega))} \leq 1/2,$
the definition of $\bV_0(\bk)$ gives us that
\begin{equation*}
\|\bV_0(\bk)\|_{L_{\infty}(0,T;L_{q_\ell}(\Omega))} 
\lesssim \|\nabla \bu \|_{L_1(0,T;L_{q_\ell}(\Omega))} 
 \lesssim \|\Jt^b \nabla \bu\|_{L_p(0,T;L_{q_\ell}(\Omega))},
 \,\,\, \forall\,\,\ell=1,2,
\end{equation*}
which, together with \eqref{es:V_1}, yields that
\begin{equation}\label{es:V_2}
\begin{aligned}
\|\bV_0(\bk)\|_{L_{\infty}(0,T;L_{\infty}(\Omega))} 
+ \sum_{\ell=1,2} \|\bV_0(\bk)\|_{L_{\infty}(0,T;H^1_{q_\ell}(\Omega))} 
\lesssim \sum_{\ell=1,2} \|\Jt^b \nabla \bu\|_{L_p(0,T;H^{1}_{q_\ell}(\Omega))}
\lesssim \CE_T.
\end{aligned}
\end{equation}
\smallbreak

Now let us derive the estimates of 
$f(\eta,\bu)= -\eta \,\dv \bu -(\gamma_1+\eta) \bV_0(\bk):\nabla \bu.$
According to the definition of $\CE_T$ and H\"older inequalities, we have
\begin{equation}\label{es:f_1_1}
\begin{aligned}
\|\Jt^b \eta \dv \bu \|_{L_p(0,T;H^{1}_{q_1/2}(\Omega))}
\lesssim & \|\eta\|_{L_\infty(0,T;L_{q_1}(\Omega))} 
\|\Jt^b \nabla \bu\|_{L_p(0,T;H^1_{q_1}(\Omega))}\\
&+\| \nabla \eta\|_{L_{\infty}(0,T;L_{q_1}(\Omega))}
\|\Jt^b \nabla \bu \|_{L_p(0,T;L_{q_1}(\Omega))}\\
\lesssim &  \|\eta\|_{L_\infty(0,T;H^1_{q_1}(\Omega))}  \CE_T,
\end{aligned}
\end{equation}
\begin{equation}\label{es:f_1_2}
\begin{aligned}
\|\Jt^b \eta \dv \bu \|_{L_p(0,T;H^{1}_{q_\ell}(\Omega))}
\lesssim & \|\eta\|_{L_\infty(0,T;L_{\infty}(\Omega))} 
\|\Jt^b \nabla \bu\|_{L_p(0,T;H^1_{q_\ell}(\Omega))}\\
&+\|\nabla \eta\|_{L_{\infty}(0,T;L_{q_\ell}(\Omega))}
\|\Jt^b \nabla \bu \|_{L_p(0,T;L_{\infty}(\Omega))}\\
\lesssim &  \big( \|\eta\|_{L_\infty(0,T;L_{\infty}(\Omega))} 
+\|\nabla \eta \|_{L_\infty(0,T;L_{q_{\ell}}(\Omega))} \big) \CE_T,
\,\,\,\ell =1,2.
\end{aligned}
\end{equation}

Next,  for the bound of $\bV_0(\bk):\nabla \bu,$ we use  H\"older inequalities and \eqref{es:V_2} to get
\begin{equation}\label{es:f_2_1}
\begin{aligned}
\|\Jt^b \bV_0(\bk):\nabla \bu \|_{L_p(0,T;H^{1}_{q_1/2}(\Omega))}
\lesssim & \|\bV_0(\bk)\|_{L_\infty(0,T;L_{q_1}(\Omega))} 
\|\Jt^b \nabla \bu\|_{L_p(0,T;H^1_{q_1}(\Omega))}\\
&+\|\nabla \bV_0(\bk) \|_{L_\infty(0,T;L_{q_1}(\Omega))}
\|\Jt^b \nabla \bu \|_{L_p(0,T;L_{q_1}(\Omega))}
\lesssim \CE_T^2,
\end{aligned}
\end{equation}
\begin{equation}\label{es:f_2_2}
\begin{aligned}
\|\Jt^b \bV_0(\bk):\nabla \bu \|_{L_p(0,T;H^{1}_{q_\ell}(\Omega))}
\lesssim & \|\bV_0(\bk)\|_{L_\infty(0,T;L_{\infty}(\Omega))} 
\|\Jt^b \nabla \bu\|_{L_p(0,T;H^1_{q_\ell}(\Omega))}\\
&+\|\nabla \bV_0(\bk) \|_{L_\infty(0,T;L_{q_\ell}(\Omega))}
\|\Jt^b \nabla \bu \|_{L_p(0,T;L_{\infty}(\Omega))}
\lesssim \CE_T^2
\end{aligned}
\end{equation}
for $\ell=1,2.$
Analogously, we have 
\begin{equation}\label{es:f_3_1}
\begin{aligned}
& \|\Jt^b \eta \bV_0(\bk):\nabla \bu \|_{L_p(0,T;H^{1}_{q_1/2}(\Omega))} \\
\lesssim &  \|\bV_0(\bk)\|_{L_\infty(0,T;L_{\infty}(\Omega))}  
 \|\eta\|_{L_\infty(0,T;L_{q_1}(\Omega))} 
\|\Jt^b \nabla \bu\|_{L_p(0,T;H^1_{q_1}(\Omega))}\\
& +\|\bV_0(\bk)\|_{L_\infty(0,T;L_{\infty}(\Omega))}  
\| \nabla \eta\|_{L_\infty(0,T;L_{q_1}(\Omega))}
\|\Jt^b\nabla \bu \|_{L_p(0,T;L_{q_1}(\Omega))}\\
& +\|\nabla \bV_0(\bk) \|_{L_\infty(0,T;L_{q_1}(\Omega))}
\|\eta\|_{L_{\infty}(0,T;L_{\infty}(\Omega))}
\|\Jt^b \nabla \bu\|_{L_p(0,T;L_{q_1}(\Omega))}\\
\lesssim &  \big( \|\eta\|_{L_\infty(0,T;H^1_{q_1}(\Omega))}
+ \|\eta\|_{L_{\infty}(0,T;L_{\infty}(\Omega))} \big)
\CE_T^2,
\end{aligned}
\end{equation}
and 
\begin{equation}\label{es:f_3_2}
\begin{aligned}
& \|\Jt^b \eta \bV_0(\bk):\nabla \bu \|_{L_p(0,T;H^{1}_{q_\ell}(\Omega))} \\
\lesssim &  \|\bV_0(\bk)\|_{L_\infty(0,T;L_{\infty}(\Omega))}  
 \|\eta\|_{L_\infty(0,T;L_{\infty}(\Omega))} 
\|\Jt^b \nabla \bu\|_{L_p(0,T;H^1_{q_\ell}(\Omega))}\\
& +\|\bV_0(\bk)\|_{L_\infty(0,T;L_{\infty}(\Omega))}  
\| \nabla \eta\|_{L_\infty(0,T;L_{q_\ell}(\Omega))}
\|\Jt^b\nabla \bu \|_{L_p(0,T;L_{\infty}(\Omega))}\\
& +\|\nabla \bV_0(\bk) \|_{L_\infty(0,T;L_{q_\ell}(\Omega))}
\|\eta\|_{L_{\infty}(0,T;L_{\infty}(\Omega))}
\|\Jt^b \nabla \bu\|_{L_p(0,T;L_{\infty}(\Omega))}\\
\lesssim & \big( \|\eta\|_{L_\infty(0,T;L_{\infty}(\Omega))} 
+\|\nabla \eta \|_{L_\infty(0,T;L_{q_{\ell}}(\Omega))} \big) \CE_T^2.
\end{aligned}
\end{equation}

Hence, summing up the bounds \eqref{es:f_1_1}, \eqref{es:f_1_2}, 
\eqref{es:f_2_1}, \eqref{es:f_2_2}, \eqref{es:f_3_1} and \eqref{es:f_3_2},  we obtain that 
\begin{equation}\label{es:f_1}
\sum_{q\in \{q_1/2,\,q_1,q_2\}}\|\Jt^b f(\eta, \bu)\|_{L_p(0,T;H^{1}_{q}(\Omega))}
\lesssim \Big( \sum_{\ell=1,2}\|\eta\|_{L_\infty(0,T;H^1_{q_\ell}(\Omega))} \Big) 
(\CE_T+\CE_T^2)+\CE_T^2.
\end{equation}
To complete the discussion of $f(\eta, \bu),$
note that $\eta(\cdot,t)=\eta(\cdot,0)+\int_0^t \pd_s \eta(\cdot,s)\,ds.$ 
Then we have
\begin{equation}\label{es:eta_infty}
\begin{aligned}
 \|\eta\|_{L_{\infty}(0,T; H^{m}_{q_\ell}(\Omega))}
\leq \|\rho_0\|_{H^{m}_{q_\ell}(\Omega)}
 +C_{p,b}\|\Jt^b \pd_t \eta\|_{L_p(0,T;H^m_{q_\ell}(\Omega))}
\end{aligned}
\end{equation}
for $\ell =1,2,$ and $m=0,1.$
Therefore, we insert \eqref{es:eta_infty} into \eqref{es:f_1} to get
\begin{equation}\label{es:f_2}
\begin{aligned}
\sum_{q\in \{q_1/2,\,q_1,q_2\}}
\|\Jt^b f(\eta, \bu)\|_{L_p(0,T;H^{1}_{q}(\Omega))}
&\lesssim \big(\|\rho_0\|_{H^1_{q_1}(\Omega)}
 +\|\rho_0\|_{H^1_{q_2}(\Omega)} \big) \CE_T+ \CE_T^2+\CE_T^3\\
& \lesssim \CI_0' + \CE_T^2+\CE_T^3,
\end{aligned}
\end{equation}
where we have used $\CI_0'\leq 1$ to obtain $\CI_0'^2 \leq \CI_0'$.

\subsection{Bound of $\bg(\eta,\bu)$}
In this subsection, we study the source term $\bg(\eta,\bu)$ defined by
\begin{equation*}
\bg(\eta,\bu)= -\eta \pd_t \bu 
+\DV\big(\CS_{\bD}(\bu) - Q(\eta)\bI\big)
+\big( \bV_0(\bk) \nabla \mid \bS_{\bu}(\bu) - \gamma_2 \eta \bI - Q(\eta)\bI\big).
\end{equation*}
More precisely, we shall prove that 
\begin{equation}\label{es:g}
\begin{aligned}
\sum_{q\in \{q_1/2,\,q_1,q_2\}}\|\Jt^b \bg(\eta, \bu)\|_{L_p(0,T;L_{q}(\Omega))}
\lesssim &\big(\|\rho_0\|_{L_{q_1}(\Omega)}
 +\|\rho_0\|_{H^1_{q_2}(\Omega)} \big)  \CE_T+ \CE_T^2+\CE_T^3\\
\lesssim  &  \CI_0' + \CE_T^2+\CE_T^3,
\end{aligned} 
\end{equation} 
since $\|\rho_0\|_{L_{q_1}(\Omega)}
 +\|\rho_0\|_{H^1_{q_2}(\Omega)} \leq \CI_0' \leq 1.$
\medskip

Firstly, by H\"older inequalities and \eqref{es:eta_infty}, we have 
\begin{equation}\label{es:g_1}
\begin{aligned}
\|\Jt^b \eta \pd_t \bu\|_{L_p(0,T;L_{q_1/2}(\Omega))}
&\lesssim \|\eta\|_{L_\infty(0,T;L_{q_1}(\Omega))}
\|\Jt^b \pd_t \bu\|_{L_p(0,T;L_{q_1}(\Omega))}
\lesssim (\|\rho_0\|_{L_{q_1}(\Omega)} + \CE_T)\CE_T,\\
\|\Jt^b \eta \pd_t \bu\|_{L_p(0,T;L_{q_\ell}(\Omega))}
&\lesssim \|\eta\|_{L_\infty(0,T;L_{\infty}(\Omega))}
\|\Jt^b \pd_t \bu\|_{L_p(0,T;L_{q_{\ell}}(\Omega))}
\lesssim (\|\rho_0\|_{H^1_{q_2}(\Omega)} + \CE_T)\CE_T,
\end{aligned}
\end{equation}
for $\ell=1,2.$
\smallbreak 

Secondly, from \eqref{es:f_2_1} and \eqref{es:f_2_2}, it follows immediately that
\begin{equation}\label{es:g_2}
\sum_{q\in \{q_1/2,\,q_1,q_2\}}
\|\Jt^b \DV \CS_{\bD} (\bu)\|_{L_p(0,T;L_{q}(\Omega))}
\lesssim \CE_T^2.
\end{equation}

Next, from the definition of $Q(\eta)$, we see easily that 
$$|\nabla Q(\eta)| \lesssim |\eta| |\nabla \eta| + |\eta|^2|\nabla \eta|.$$
Then, we obtain from \eqref{cdt:etau} and \eqref{es:eta_infty} that
\begin{equation}\label{es:g_3_1}
\begin{aligned}
\|\Jt^b \nabla Q(\eta)\|_{L_p(0,T;L_{q_1/2}(\Omega))}
\lesssim & \big(1+\|\eta\|_{L_{\infty}(0,T;L_{\infty}(\Omega))}\big)
\|\eta\|_{L_{\infty}(0,T;L_{q_1}(\Omega))}
\|\Jt^b \nabla \eta\|_{L_{p}(0,T;L_{q_1}(\Omega))}\\
\lesssim & (\|\rho_0\|_{L_{q_1}(\Omega)}+\CE_T) \CE_T,\\
\|\Jt^b \nabla Q(\eta)\|_{L_p(0,T;L_{q_\ell}(\Omega))}
\lesssim & \big(1+\|\eta\|_{L_{\infty}(0,T;L_{\infty}(\Omega))}\big)
\|\eta\|_{L_{\infty}(0,T;L_{\infty}(\Omega))}
\|\Jt^b \nabla \eta\|_{L_{p}(0,T;L_{q_{\ell}}(\Omega))} \\
\lesssim & (\|\rho_0\|_{H^1_{q_2}(\Omega)}+\CE_T) \CE_T,
\,\,\,\forall\,\,\ell=1,2,
\end{aligned}
\end{equation}
which give us 
\begin{equation}\label{es:g_3}
\sum_{q\in \{q_1/2,\,q_1,q_2\}}
\|\Jt^b \nabla Q(\eta)\|_{L_p(0,T;L_{q}(\Omega))}
\lesssim \big(\|\rho_0\|_{L_{q_1}(\Omega)}
 +\|\rho_0\|_{H^1_{q_2}(\Omega)} \big)  \CE_T+\CE_T^2.
\end{equation}

As $\bS_{\bu}(\bu)=\bS(\bu)+\CS_{\bD}(\bu)$, 
by \eqref{es:V_2} and \eqref{es:g_2} we have 
\begin{equation}\label{es:g_4_1}
\begin{aligned}
&\big\|\Jt^b\big( \bV_0(\bk) \nabla \mid \bS_{\bu}(\bu) -\gamma_2 \eta\big)\big\|_{L_p(0,T;L_{q_1/2}(\Omega))} \\
\lesssim  & \|\bV_0(\bk)\|_{L_\infty(0,T;L_{q_1}(\Omega))} 
\big(   \|\Jt^b  (\nabla^2\bu,\nabla \eta)\|_{L_p(0,T;L_{q_1}(\Omega))} 
+ \|\Jt^b \nabla \CS_{\bD} (\bu)\|_{L_p(0,T;L_{q_1}(\Omega))} \big) \\
\lesssim &  \CE_T^2+\CE_T^3,\\
&\big\|\Jt^b\big( \bV_0(\bk) \nabla \mid \bS_{\bu}(\bu)-\gamma_2 \eta\big)\big\|_{L_p(0,T;L_{q_\ell}(\Omega))} \\
\lesssim  & \|\bV_0(\bk)\|_{L_\infty(0,T;L_{\infty}(\Omega))} 
\big(   \|\Jt^b (\nabla^2 \bu,\nabla \eta)\|_{L_p(0,T;L_{q_\ell}(\Omega))} 
+ \|\Jt^b \nabla \CS_{\bD} (\bu)\|_{L_p(0,T;L_{q_\ell}(\Omega))} \big) \\
\lesssim &  \CE_T^2+\CE_T^3.
\end{aligned}
\end{equation}
And, by \eqref{es:V_2} and \eqref{es:g_3_1} we have 
\begin{equation}\label{es:g_4_2}
\begin{aligned}
\big\|\Jt^b\big( \bV_0(\bk) \nabla \mid Q(\eta)\bI\big)\big\|_{L_p(0,T;L_{q_1/2}(\Omega))} \lesssim  & \|\bV_0(\bk)\|_{L_\infty(0,T;L_{q_1}(\Omega))} 
\|\Jt^b \nabla Q (\eta)\|_{L_p(0,T;L_{q_1}(\Omega))} \\
\lesssim & \|\rho_0\|_{H^1_{q_2}(\Omega)}  \CE_T^2+\CE_T^3,\\
\big\|\Jt^b\big( \bV_0(\bk) \nabla \mid Q(\eta)\bI\big)\big\|_{L_p(0,T;L_{q_\ell}(\Omega))} \lesssim  & \|\bV_0(\bk)\|_{L_\infty(0,T;L_{\infty}(\Omega))} 
\|\Jt^b \nabla Q (\eta)\|_{L_p(0,T;L_{q_\ell}(\Omega))} \\
\lesssim & \|\rho_0\|_{H^1_{q_2}(\Omega)}  \CE_T^2+\CE_T^3.
\end{aligned}
\end{equation}
Combing \eqref{es:g_4_1} and \eqref{es:g_4_2} yields that 
\begin{equation}\label{es:g_4}
\sum_{q\in \{q_1/2,\,q_1,q_2\}}
\big\|\Jt^b \big( \bV_0(\bk) \nabla \mid \bS_{\bu}(\bu) - \gamma_2 \eta \bI - Q(\eta)\bI\big)\big\|_{L_p(0,T;L_{q}(\Omega))}
\lesssim  \CE_T^2+\CE_T^3.
\end{equation}
Therefore, we can conclude the claim \eqref{es:g}
from \eqref{es:g_1}, \eqref{es:g_2}, \eqref{es:g_3} and \eqref{es:g_4}.

\subsection{Bound of the boundary term}

In order to estimate the boundary term, first of all we introduce the following
extension formula of functions defined on $(0, T)$.  Let $f$ be a function defined
on $(0, T)$ such that $f|_{t=0}=0$. 
We set 
$$e[f](\cdot, t) = \begin{cases}
 f(\cdot, t) &\quad t \in (0, T), \\
f(\cdot, 2T-t) &\quad t \in (T, 2T), \\
0 &\quad t \not\in (0, 2T).
\end{cases}$$
From the definition, we have
\begin{equation}\label{eq:Dt_e}
\pd_te[f](\cdot, t) = \begin{cases} 
(\pd_tf)(\cdot, t) &\quad t \in (0, T), \\
-(\pd_tf)(\cdot, 2T-t)  &\quad t \in (T, 2T), \\
0 &\quad t \not\in (0, 2T).
\end{cases}
\end{equation}

Since $\eta|_{t=0} \not=0$, we introduce another extension formula of $\eta$.
Let $\eta_L \in H^1_p(\BR, H^1_{q_1}(\Omega) \cap H^1_{q_2}(\Omega))$ be a 
function such that $\eta_L|_{t=0} = \rho_0$ and 
\begin{equation}\label{es:eta_L}
\|e^{\gamma |t|}\eta_L\|_{H^1_p(\BR, H^1_{q_\ell}(\Omega))}
\leq C\|\rho_0\|_{H^1_{q_\ell}(\Omega)} \quad (\ell=1,2)
\end{equation}
for some constants $\gamma,C>0.$
Such $\eta_L$ can be constructed as follows. Let $\tilde \rho_0$ be an extension
of $\rho_0$ to $\BR^3$ such that $\tilde \rho_0|_{\Omega} = \rho_0$ and 
$$\|\tilde\rho_0\|_{H^1_{q_{\ell}}(\BR^3)} 
\leq C\|\rho_0\|_{H^1_{q_\ell}(\Omega)},\quad 
 \forall\,\,\, \ell=1,2. $$
Then, let $(\rho_L,\bv_L)$ be a solution of the equations for some $\lambda_0>0$:
\begin{equation*} 
	\left\{ \begin{aligned}
&\pd_t \rho_L + \lambda_0\rho_L + \gamma_1 \dv v_L =0
&&\quad&\text{in} &\quad \BR^3\times \BR_+, \\
&\gamma_1 \pd_t\bv_L + \lambda_0\bv_L- \mu\Delta \bv_L 
-\nu\nabla\dv\bv_L
+ \gamma_2 \nabla \rho_L=0
&&\quad&\text{in}& \quad \BR^3\times \BR_+, \\
&(\rho_L, \bv_L)|_{t=0} = (\tilde\rho_0, 0)
&&\quad&\text{in}& \quad \BR^3
\end{aligned}
\right.
\end{equation*}
And then, set $\eta_{L}(t) = \rho_L(|t|)$.
\smallbreak 

Now, we define the extension $\Ext[\bh(\eta, \bu)]$ of $\bh(\eta, \bu)$ as follows. Let 
\begin{align*}
\CS_\bD^e &= \mu e[\Jt^b\CD_\bD(\bu)]\bn_\Gamma
+ (\nu-\mu)e[\Jt^b\bV_0(\bk):\nabla\bu]\bn_\Gamma, \\
Q^e &= \Jt^bQ(\eta_L) 
+ e\big[ \Jt^b\big( Q(\eta) - Q(\eta_L) \big)\big]\bn_\Gamma,  \\
\CS_\bu^e & = e[\Jt^b\bS_\bu(\bu)\bV_0(\bk)]\bv_\Gamma, \\
R^e & =e[\Jt^b(\gamma_2\eta+Q(\eta))\bV_0(\bk)]\bn_\Gamma.
\end{align*}
And then, set 
$\Ext[\bh(\eta, \bu)] 
=-\CS_\bD^e + Q^e - \CS_\bu^e+ R^e.$  
Obviously, 
$$\Jt^{-b}\Ext[\bh(\eta, \bu)] = \bh(\eta, \bu) 
\,\,\, \text{for}\,\,\,  t \in (0, T). $$
Note that 
\begin{equation}\label{es:q2_1}
\begin{aligned}
& \sum_{q\in \{q_1/2,\, q_1, q_2\}} 
  \Big( \|\kappa \Ext[\bh(\eta, \bu)] \|_{L_p(\BR; H^1_q(\Omega))}
+ \|\kappa \Ext[\bh(\eta, \bu)] \|_{H^{1/2}_p(\BR; L_q(\Omega))}\Big)\\
& \lesssim \|\Ext[\bh(\eta, \bu)] \|_{L_p(\BR; H^1_{q_2}(\Omega))}
+ \|\Ext[\bh(\eta, \bu)] \|_{H^{1/2}_p(\BR; L_{q_2}(\Omega))}.
\end{aligned}
\end{equation}
In order to estimate the right-hand side of \eqref{es:q2_1}, we use the following technical lemma whose proof is given in Sect. 3.4 below. 
\begin{lem}\label{lem:half} 
Assume that $\Omega$ is a uniform $C^2$ domain in $\BR^N.$ Let 
$1 < p < \infty$ and $N < q < \infty$. Let $f$ and $g$ be functions defined for $t\in (0,T)$ fulfilling $f|_{t=0}=0.$ Then 
\begin{align*}
&\quad \|e[f\nabla g]\|_{H^{1/2}_p(\BR; L_q(\Omega))}
+ \|e[f\nabla g]\|_{L_p(\BR; H^1_q(\Omega))} \\
&\leq C\Big(
\|\pd_tf\|_{L_p(0, T; H^1_q(\Omega))}
\| g\|_{L_\infty(0, T; L_q(\Omega))}
+ \|f\|_{L_\infty(0, T; H^1_q(\Omega))}
\big(\|\pd_tg\|_{L_p(0, T; L_q(\Omega))}
+ \|g\|_{L_p(0,T; H^2_q(\Omega))} \big)\Big).
\end{align*}
\end{lem}

By Lemma \ref{lem:es_V} and \eqref{es:V_2}, we have 
\begin{equation}\label{es:V_3}
\|\pd_t \bV_0 (\bk)\|_{L_p(0, T; H^1_{q_2}(\Omega))}
+\|\bV_0 (\bk)\|_{L_\infty(0, T; H^1_{q_2}(\Omega))} 
\lesssim \| \Jt^b \bu \|_{L_p(0, T; H^1_{q_2}(\Omega))} 
\lesssim \CE_{T,q_2},
\end{equation}
for $\CE_{T,q_2} =  \|\Jt^b(\eta, \bu)\|_{L_\infty(0, T; L_{q_2}(\Omega))}  
+\CM^{b}_{p,q_{2}}(T;\eta,\bu).$ In particular, \eqref{es:V_3} also yields that 
\begin{equation}\label{es:V_4}
\|\pd_t [\bV_0 (\bk) \otimes\bV_0 (\bk)]\|_{L_p(0, T; H^1_{q_2}(\Omega))}
+\|\bV_0 (\bk) \otimes\bV_0 (\bk) \|_{L_\infty(0, T; H^1_{q_2}(\Omega))} 
\lesssim \CE_{T,q_2}^2,
\end{equation}
as $H^{1}_{q_2}(\Omega)$ is an algebra for $q_2>3.$
Here, $\bA \otimes \bB=[A_{ij} B_{k\ell}]$ for any $3\times 3$ matrices $\bA,\bB.$
Then , 
Lemma \ref{lem:half}, \eqref{es:V_3}, \eqref{es:V_4} and \eqref{es:Jt_1} yields that 
\begin{equation}\label{es:q2_2}
\begin{aligned}
 \|\CS_\bD^e\|_{L_p(\BR; H^1_{q_2}(\Omega))}
+ \|\CS_\bD^e\|_{H^{1/2}_p(\BR; L_{q_2}(\Omega))}
& \lesssim \CE_{T,q_2}^2,\\
 \|\CS_\bu^e\|_{L_p(\BR; H^1_{q_2}(\Omega))}
+ \|\CS_\bu^e\|_{H^{1/2}_p(\BR; L_{q_2}(\Omega))}
& \lesssim \CE_{T,q_2}^3.
\end{aligned}
\end{equation}

For the estimate of $Q^e$ and $R^e$, we use the continuous embedding
$$H^1_p\big(\BR; L_q(\Omega)\big) \hookrightarrow 
H^{1/2}_p \big(\BR; L_q(\Omega)\big).$$ 
According to \eqref{es:eta_L}, we have 
\begin{equation*}
\|\eta_{L}\|_{L_\infty(\BR;L_{\infty}(\Omega))} \lesssim \|\rho_0\|_{H^1_{q_2}(\Omega)},
\end{equation*}
which yields that 
\begin{equation}\label{es:Q_1}
\| \Jt^bQ(\eta_L)\|_{L_p(\BR; H^1_{q_2}(\Omega))} 
+\|\Jt^bQ(\eta_L)\|_{H^1_p(\BR;L_{q_2}(\Omega))} 
\lesssim \|\rho_0\|_{H^1_{q_2}(\Omega)}^2.
\end{equation}
Notice that $\|\eta\|_{L_\infty(0, T; L_{q_2})} \leq \|\Jt^b\eta\|_{L_\infty(0, T; L_{q_2})}
\leq \CE_{T, q_2}$. And then, 
\begin{equation}\label{es:Q_2}
\begin{aligned}
&\| \Jt^bQ(\eta)\|_{L_p(0,T; H^1_{q_2}(\Omega))} 
+\|\Jt^b \pd_t Q(\eta)\|_{L_p(0,T;L_{q_2}(\Omega))} \\
\lesssim &  \|\Jt^b\eta\|_{L_p(0, T; H^1_{q_2}(\Omega))}
\|\eta\|_{L_\infty(0, T; H^1_{q_2}(\Omega))} 
+\|\Jt^b\pd_t\eta\|_{L_p(0,T;H^1_{q_2}(\Omega))}
\|\eta\|_{L_\infty(0, T; L_{q_2}(\Omega))} \\
\lesssim &  \CE_{T, q_2} \big(\|\rho_0\|_{H^1_{q_2}(\Omega)}
+\|\Jt^b\pd_t\eta\|_{L_p(0,T;H^1_{q_2}(\Omega))} \big)
+\CE_{T, q_2}^2\lesssim \|\rho_0\|_{H^1_{q_2}(\Omega)}^2+\CE_{T, q_2}^2.
\end{aligned}
\end{equation}
Thus we use \eqref{es:Q_1} and \eqref{es:Q_2} to obtain
\begin{equation}\label{es:q2_3}
\begin{aligned}
&\|Q^e\|_{L_p(\BR; H^1_{q_2}(\Omega))} 
+\|Q^e\|_{H^{1/2}_p(\BR;L_{q_2}(\Omega))} \\
\lesssim & \|Q^e\|_{L_p(\BR; H^1_{q_2}(\Omega))}
+\|\pd_t Q^e\|_{L_p(\BR; L_{q_2}(\Omega))} \\
\lesssim & \CI_0^2
+ \big\|\Jt^b\big(Q(\eta) - Q(\eta_L)\big)\|_{L_p(0, T;H^1_{q_2}(\Omega))} 
+ \big\|\Jt^b \pd_t \big(Q(\eta) - Q(\eta_L)\big)\|_{L_p(0, T; L_{q_2}(\Omega))}\\
\lesssim  &  \|\rho_0\|_{H^1_{q_2}(\Omega)}^2 + \CE_{T,q_2}^2.
\end{aligned}
\end{equation}

Analogously, by \eqref{es:V_3} and \eqref{es:Q_2} we have 
\begin{equation}\label{es:q2_4}
\begin{aligned}
&\|R^e\|_{L_p(\BR; H^1_{q_2}(\Omega))}
+\|R^e\|_{H^{1/2}_p(\BR;L_{q_2}(\Omega))}  \\
 \lesssim & \|R^e\|_{L_p(\BR; H^1_{q_2}(\Omega))}
+\|\pd_t R^e\|_{L_p(\BR; L_{q_2}(\Omega))} \\
\lesssim & \big\|\Jt^b \big(\gamma_2\eta+Q(\eta)\big)
\bV_0(\bk)\big\|_{L_p(0,T; H^1_{q_2}(\Omega))} 
+\big\|\pd_t\big( \Jt^b(\gamma_2\eta+Q(\eta))
\bV_0(\bk) \big) \big\|_{L_p(0, T; L_{q_2}(\Omega))}\\
\lesssim & \CE_{T,q_2}^2+\big\|\Jt^b Q(\eta)
\bV_0(\bk)\big\|_{L_p(0,T; H^1_{q_2}(\Omega))} 
+\big\| \Jt^b \pd_t\big(Q(\eta)\bV_0(\bk) \big) \big\|_{L_p(0, T; L_{q_2}(\Omega))}\\
\lesssim &  \CE_{T,q_2}^2+\big( \|\rho_0\|_{H^1_{q_2}(\Omega)}^2+\CE_{T, q_2}^2 \big) \CE_{T,q_2} +\|\Jt^bQ(\eta)\|_{L_\infty(0, T; L_{q_2}(\Omega))}
\|\pd_t \bV_0(\bk)\|_{L_p(0, T; L_{\infty}(\Omega))}\\
\lesssim &\|\rho_0\|_{H^1_{q_2}} + \CE_{T,q_2}^2+\CE_{T,q_2}^3.
\end{aligned}
\end{equation}
Here, we have used the facts that $ \|\rho_0\|_{H^{1}_{q_2}(\Omega)} \leq \CI_0'\leq 1$ and that
$$ \|\eta\|_{L_{\infty}(0,T; H^{1}_{q_2}(\Omega))}
\lesssim \|\rho_0\|_{H^{1}_{q_2}(\Omega)}
+\|\Jt^b \pd_t \eta\|_{L_p(0,T;H^1_{q_2}(\Omega))} 
 \lesssim \|\rho_0\|_{H^{1}_{q_2}(\Omega)} + \CE_{T, q_2}.$$ 

At last, combining all the bounds \eqref{es:q2_1},  \eqref{es:q2_2} \eqref{es:q2_3} and \eqref{es:q2_4} yields that 
\begin{equation}\label{es:q2}
\begin{aligned}
& \sum_{q\in \{q_1/2,\, q_1, q_2\}} 
  \Big( \|\kappa \Ext[\bh(\eta, \bu)] \|_{L_p(\BR; H^1_q(\Omega))}
+ \|\kappa \Ext[\bh(\eta, \bu)] \|_{H^{1/2}_p(\BR; L_q(\Omega))}\Big)\\
& \lesssim \|\rho_0\|_{H^1_{q_2}} + \CE_{T,q_2}^2+\CE_{T,q_2}^3.
\end{aligned}
\end{equation}
because the support of $\kappa$ is compact and $q_1/2 < q_1 < q_2$.
\subsection{Proof of Lemma \ref{lem:half}}
We end up with the proof of Lemma \ref{lem:half}. 
To this end, we introduce the dual space of 
$H^1_{q'}(\BR^N)$ ($1<q<\infty$), that is,
$$H^{-1}_q(\BR^N) = \{f \in \CS'(\BR^N) \cap L_{1, {\rm loc}}(\BR^N) \mid
\|f\|_{H^{-1}_q(\BR^N)} < \infty\}$$ with the norm given by
$$\|f\|_{H^{-1}_q(\BR^N)} =\big\|\CF^{-1}_{\xi} 
\big[(1+|\xi|^2)^{-1/2}\CF_x[f](\xi) \big]\big\|_{L_q(\BR^N)}.$$
Here, the symbols $\CF_{x}$ and $\CF^{-1}_{\xi}$ denote the Fourier transformation and its inverse in $\BR^N.$

\begin{proof}[Proof of Lemma \ref{lem:half}]

By extending functions defined on $\Omega$ to $\BR^N$, we may assume that 
$\Omega=\BR^N.$ Then we know that 
$$\|h\|_{H^{1/2}_p(\BR; L_q(\BR^N))} \lesssim 
\|\pd_th\|_{H^1_p(\BR;H^{-1}_q(\BR^N))}
+ \|h\|_{L_p(\BR; H^1_q(\BR^N))},$$
which together with Sobolev inequalities yields that
\begin{equation}\label{es:fDg_1}
\begin{aligned}
 \|e[f\nabla g]\|_{H^{1/2}_p(\BR; L_q(\BR^N))} 
& \lesssim \|\pd_t e[f\nabla g]\|_{L_p(\BR; H^{-1}_q(\BR^N))}
+ \|e[f\nabla g]\|_{L_p(\BR; H^1_q(\BR^N))}\\
&\lesssim \|\pd_t e[f\nabla g]\|_{L_p(\BR; H^{-1}_q(\BR^N))}
+  \|f\|_{L_\infty(0,T; H^1_q(\BR^N))}
\|\nabla g\|_{L_p(0,T; H^1_q(\BR^N))}.
\end{aligned}
\end{equation}

Since $f\nabla g|_{t=0},$ we infer from \eqref{eq:Dt_e} that 
$$\pd_t\, e [f\nabla g]  (\cdot, t)
=\begin{cases}
 \big( \pd_t (f \nabla g) \big)  (\cdot, t), &t\in (0, T) ,\\
 - \big( \pd_t (f \nabla g) \big)  (\cdot, 2T-t), & t\in (T, 2T), \\
0 & t \not\in (0, 2T).
\end{cases}$$

For any $\varphi \in C^\infty_0(\BR^N),$ $0<t<T$ and $k=1,\dots,N,$ 
integration by parts and H\"older inequalities imply that 
\begin{align*}
\Big|\big(\pd_t(f\pd_k g), \varphi \big)_{H^{-1}_q(\BR^N)
\times H^1_{q'}(\BR^N)}  \Big| 
\leq & \Big| \int_{\BR^N} (\pd_t f \varphi) \, \pd_k g \,dx \Big|  
+ \Big| \int_{\BR^N} \pd_k \pd_tg f\varphi \,dx \Big|   \\
\leq &  \Big| \int_{\BR^N} \pd_k (\pd_t f \varphi) \,  g \,dx \Big|  
+ \Big| \int_{\BR^N}  \pd_tg \pd_k  (f\varphi) \,dx \Big|  \\
\leq & \Big( \|\pd_t \pd_k f\|_{L_q(\BR^N)}
\| g\|_{L_q(\BR^N)}
+\|\pd_k f\|_{L_q(\BR^N)}
\| \pd_t  g\|_{L_q(\BR^N)} \Big)
\|\varphi\|_{L_s(\BR^N)}\\
&+\Big( \|\pd_t f\|_{L_\infty(\BR^N)}
\| g\|_{L_q(\BR^N)}
+\| f\|_{L_\infty(\BR^N)}
\| \pd_t g\|_{L_q(\BR^N)}
\Big)\| \pd_k \varphi\|_{L_{q'}(\BR^N)}.
\end{align*}
with $s$ satisfying $2/q + 1/s=1.$  
Since 
$$N(1/q'-1/s)
= N\big( 1-1/q-(1-2/q) \big) = N/q < 1,$$
we have the Sobolev's inequality
$$\|\varphi\|_{L_s(\BR^N)} \leq C\|\varphi\|_{H^1_{q'}(\BR^N)}. $$
Combining these estimates yields that 
\begin{equation*}
\begin{aligned}
\|\pd_t(f\nabla g)(\cdot, t)\|_{H^{-1}_q(\BR^N)}
\leq & \sup_{k=1,\dots,N}\sup_{\varphi \in \CB_1}\Big|\big(\pd_t(f\pd_k g), \varphi \big)_{H^{-1}_q(\BR^N)
\times H^1_{q'}(\BR^N)}  \Big| \\
\lesssim & \|\pd_t f\|_{H^1_q(\BR^N)}\|g\|_{L_q(\BR^N)} 
+ \|f\|_{H^1_q(\BR^N)}\|\pd_tg\|_{L_q(\BR^N)}.
\end{aligned}
\end{equation*}
for $\CB_1=\{\varphi \in H^1_{q'}(\BR^N) : 
 \|\varphi\|_{H^1_{q'}(\BR^N)} \leq 1\}$ 
and $0<t<T.$ 
Therefore,
\begin{equation}\label{es:fDg_2}
\begin{aligned}
\|\pd_t\,e [f\nabla g]\|_{L_p(\BR; H^{-1}_q(\BR^N))}
\lesssim & \|\pd_tf\|_{L_p(0,T; H^1_q(\BR^N))}
\|g\|_{L_\infty(0,T; L_q(\BR^N))} \\
& \quad +\|f\|_{L_\infty(0,T; H^1_q(\BR^N))}
 \|\pd_tg\|_{L_p(0,T;L_q(\BR^N))}.
\end{aligned}
\end{equation}
Then \eqref{es:fDg_1} and \eqref{es:fDg_2} give us 
\begin{equation*}
\begin{aligned}
\|e[f\nabla g]\|_{H^{1/2}_p(\BR; L_q(\BR^N))} 
\lesssim & \|\pd_tf\|_{L_p(0,T; H^1_q(\BR^N))}
\|g\|_{L_\infty(0,T; L_q(\BR^N))} \\
& \quad +\|f\|_{L_\infty(0,T; H^1_q(\BR^N))}
\Big( \|\pd_tg\|_{L_p(0,T;L_q(\BR^N))} 
+\|\nabla g\|_{L_p(0,T; H^1_q(\BR^N))} \Big)
\end{aligned}
\end{equation*}
for $q>N$. This completes the proof. 
\end{proof}

\section{Some auxiliary system}
In this short section, we consider the following auxiliary system:
\begin{equation} \label{sft**}
	\left\{ \begin{aligned}
&\pd_t\rho^2 + \lambda_0 \rho^2+ \gamma_1 \, \dv\bv^2= (\lambda_0\rho^1)_0
&&\quad&\text{in} &\quad \Omega \times \BR, \\
&\gamma_1 \pd_t\bv^2 +\lambda_0 \bv^2
-\DV\big(\bS( \bv^2 ) - \gamma_2 \rho^2 \bI\big)= (\lambda_0 \bv^1)_0
&&\quad&\text{in}& \quad \Omega  \times \BR, \\
&\big(\bS(\bv^2) - \gamma_2 \rho^2 \bI\big)\bn_{\Gamma} =0
&&\quad&\text{on}& \quad \Gamma \times \BR.
\end{aligned}
\right.
\end{equation}

According to  Theorem \ref{thm:linear} and Theorem \ref{thm:nonlinear}, the system \eqref{sft**} admits a unique solution $(\rho^2,\bv^2)$ with
\begin{align*}
\rho^2 \in H^1_p\big(\BR; H^1_q(\Omega)\big),
\quad \bv^2 \in H^1_p\big(\BR; L_q(\Omega)^3\big) 
\cap L_p\big(\BR; H^2_q(\Omega)^3\big)
\end{align*}
for $q=q_1/2, q_1$ and $q_2.$ Here, the indices $b,$ $q_1$ and $q_2$ satisfy the assumptions in Theorem \ref{thm:nonlinear}. 
Moreover, the solution $(\rho^2,\bv^2)$ 
satisfies the estimates
\begin{equation}\label{es:rhov_2}
\sum_{q\in \{q_1/2,\, q_1, q_2\}} \CM^b_{p,q}(T; \rho^2, \bv^2) 
\lesssim \CI_0'+ \CE^2_T + \CE^3_T
\end{equation}
for  $\CI_0' = \sum_{\ell=1, 2} \|\rho_0\|_{H^1_{q_\ell}(\Omega)},$ and $(\rho^2, \bv^2) = (0, 0)$ for $t \leq 0.$ Thus we obtain that 
\begin{equation}\label{es:rhov_12_1}
\sum_{\fa\in\{1,2\}}\sum_{q\in \{q_1/2,\, q_1, q_2\}} \CM^b_{p,q}(T; \rho^{\fa}, \bv^{\fa}) 
\lesssim \CI_0'+ \CE^2_T + \CE^3_T.
\end{equation}
which together with Sobolev inequalities imply that
\begin{equation}\label{es:rhov_12_2}
\sum_{\fa\in\{1,2\}}\|\Jt^{b}(\rho^{\fa},\bv^{\fa})\|_{L_{\infty}(0,T;H^{1,0}_{q_1}(\Omega))}
\lesssim \CI_0'+ \CE^2_T + \CE^3_T.
\end{equation}
Therefore, assuming that $b\geq 3/(2q_1),$ we infer from  
\eqref{es:rhov_12_1} and \eqref{es:rhov_12_2} that
\begin{equation}\label{es:rhov_12}
\CE_T (\rho^1,\bv^1)+\CE_T (\rho^2,\bv^2)
 \lesssim   \CI_0' +\CE_T^2 +\CE_T^3.
\end{equation}

\section{Decay estimates related to the modified initial data}
\label{sec:mi}

 The compensation $(\theta,\bw)$ is constructed as the solution of the equations: 
\begin{equation} \label{sft***}
	\left\{ \begin{aligned}
&\pd_t\theta + \gamma_1 \, \dv\bw= \lambda_0\rho^2
&&\quad&\text{in} &\quad \Omega \times (0, T), \\
&\gamma_1 \pd_t\bw 
-\DV\big(\bS( \bw ) - \gamma_2 \theta \bI\big)= \lambda_0\bv^2
&&\quad&\text{in}& \quad \Omega  \times  (0, T), \\
&\big(\bS(\bw) - \gamma_2 \theta \bI\big)\bn_{\Gamma} =0
&&\quad&\text{on}& \quad \Gamma \times (0, T), \\
& (\theta, \bw)|_{t=0} = (\rho_0, \bv_0)-(\rho^1, \bv^1)|_{t=0}&&\quad&\text{in}&\quad\Omega.
\end{aligned}\right.
\end{equation}
From the definition of $\bh(\eta,\bu)$ and  \eqref{cdt:initial},  $(\rho^1, \bv^1)|_{t=0}$ satisfies the compatibility conditions 
$$\big(\bS(\bv^1|_{t=0}) - \gamma_2 \rho^1|_{t=0} \bI\big)\bn_{\Gamma} =Q(\rho_0)\bn_{\Gamma} = \big(\bS(\bv_0) - \gamma_2 \rho_0 \bI\big)\bn_{\Gamma}.$$
In other words, if we denote $(\theta_0,\bw_0)= (\rho_0, \bv_0)-(\rho^1, \bv^1)|_{t=0},$ then 
it holds 
\begin{equation*}
\big(\bS(\bw_0) - \gamma_2 \theta_0 \bI\big)\bn_{\Gamma} =0.
\end{equation*}
Moreover, the following estimate follows from the trace method and Theorem \ref{thm:nonlinear}
\begin{equation}\label{es:I_0_2}
\CI_0''=\sum_{q\in \{q_1,q_2\}} 
\Big( \|\theta_0\|_{H^1_q(\Omega)} 
+ \|\bw_0\|_{B^{2-2/p}_{q,p}(\Omega)} \Big)
+ \|(\theta_0, \bw_0)\|_{L_{q_1/2}(\Omega)}
\lesssim \CI_0 +\CE_T^2 +\CE_T^3.
\end{equation}
In this section, we only verify the decay estimates from $(\theta_0,\bw_0),$
while the contribution of  $(\lambda_0\rho^2,\lambda_0\bv^2)$ on the right-hand side of \eqref{sft***} will be discussed in the next section.
For convenience, we adopt the notation 
\begin{equation}\label{def:Hmn}
H^{m,n}_p(\Omega) = \{(\rho,\bv): \rho \in H^m_p(\Omega), \bv \in H^n_p(\Omega)^3\},
\quad
\|(\rho, \bv)\|_{H^{m,n}_p(\Omega)}
= \|\rho\|_{H^m_p(\Omega)} + \|\bv\|_{H^n_p(\Omega)}
\end{equation}
for $1\leq p\leq \infty$ and $m,n=0,1,2.$

\subsection{Some property of the semigroup $T(t)$}
To study the system \eqref{sft***}, we introduce the semigroup $T(t)$ associated to the following homogeneous problem
\begin{equation}\label{eq:T} 
	\left\{ \begin{aligned}
&\pd_t \ol\rho +  \gamma_1 \, \dv \ol\bv= 0
&&\quad&\text{in} &\quad \Omega \times (0,T), \\
&\gamma_1 \pd_t\ol\bv 
-\DV\big(\bS( \ol\bv ) - \gamma_2 \ol\rho \bI\big)= 0
&&\quad&\text{in}& \quad \Omega  \times (0,T), \\
&\big(\bS(\ol\bv) - \gamma_2 \ol\rho \bI\big)\bn_{\Gamma} = 0
&&\quad&\text{on}& \quad \Gamma \times (0,T), \\
&(\ol\rho, \ol\bv )|_{t=0} =(\ol\rho_0,\ol\bv_0)
&&\quad&\text{in}& \quad \Omega.
\end{aligned}
\right.
\end{equation}
We recall the semigroup setting of the equations \eqref{eq:T}.
For $1 < p < \infty$, let 
\begin{align*}
D(A) & = \big\{(\rho, \bv) \in H^{1,2}_p(\Omega) \mid 
\big(\bS(\bv) - \gamma_2\rho \bI\big)\bn_{\Gamma} = 0
\quad\text{on $\Gamma$} \,\big\}, \\
A(\rho, \bv) & =\Big(-\gamma_1\dv\bv, \gamma_1^{-1}\DV(\bS(\bv)-
\gamma_2\rho\bI\big) \Big)
\,\,\text{for $(\rho, \bv) \in D(A)$}.
\end{align*}
Then, equations \eqref{eq:T} are written as
$$\pd_t(\ol\rho, \ol\bv) - A(\ol\rho, \ol\bv) = (0, 0) \,\,\text{for
$t \in (0, T)$}, \quad(\ol\rho, \ol\bv)|_{t=0} = (\ol\rho_0, \ol\bv_0).
$$
According to the results in \cite{EvBS2014,SSZ2020}, we know that 
the operator $A$ generates a continuous analytic semigroup
$\{T(t)\}_{t\geq 0}$ on $H^{1,0}_p(\Omega)$, and therefore,
 thanks to the Duhamel principle, the solution of \eqref{sft***} is given by 
$(\theta, \bw)=(\theta^1, \bw^1) +(\theta^2, \bw^2)$ with 
\begin{equation*}
(\theta^1,\bw^1)= T(t) (\theta_0,\bw_0), \quad 
(\theta^2,\bw^2)=\lambda_0\int^t_0T(t-s)(\rho^2, \bv^2)(\cdot, s)\,ds.
\end{equation*}

Now, we review some properties of $T(t)$ used in this paper. 
From the fundamental property of the continuous analytic semigroup
we have
\begin{equation}\label{es:sg_1}
\begin{aligned}
\|T(t)(\ol\rho_0, \ol\bv_0)\|_{H^{1,0}_p(\Omega)}
&\leq C \|(\ol\rho_0, \ol\bv_0)\|_{H^{1,0}_p(\Omega)}
\quad\text{for $(\ol\rho_0, \ol\bv_0) \in H^{1,0}_p(\Omega)$}, \\
\|T(t)(\ol\rho_0, \ol\bv_0)\|_{H^{1,2}_p(\Omega)}
& \leq C\|(\ol\rho_0, \ol\bv_0)\|_{H^{1,2}_p(\Omega)}
\quad\text{for $(\ol\rho_0, \ol\bv_0) \in D(A)$}
\end{aligned}
\end{equation}
for any $0 < t \leq 2$.
To establish the long time estimates, 
we also need the following results proved in \cite{ShiZ2020a}:
\begin{equation}\label{es:sg_2}
\begin{aligned}
\|T(t)(\ol\rho_0, \ol\bv_0)\|_{L_p(\Omega)}
& \leq Ct^{-(3/q-3/p)/2}\vertiii{(\ol\rho_0,\ol\bv_0)}_{p,q}
&\,\,\,(t \geq 1), \\
\|\nabla T(t)(\ol\rho_0, \ol\bv_0)\|_{L_p(\Omega)}
& \leq Ct^{-\sigma(p,q)}\vertiii{(\ol\rho_0,\ol\bv_0)}_{p,q}
&\,\,\,(t \geq 1), \\
\|\nabla^2 \CP_v T(t)(\ol\rho_0, \ol\bv_0)\|_{L_p(\Omega)}
& \leq Ct^{-3/(2q)}\vertiii{(\ol\rho_0,\ol\bv_0)}_{p,q}
&\,\,\,(t \geq 1), 
\end{aligned}
\end{equation}
with $1\leq q \leq 2\leq p <\infty,$ 
$\vertiii{(\ol\rho_0,\ol\bv_0)}_{p,q}=\|(\ol\rho_0,\ol\bv_0)\|_{L_q(\Omega)}
+ \|(\ol\rho_0,\ol\bv_0)\|_{H^{1,0}_p(\Omega)},$ $\ol v=\CP_vT(t)(\ol\rho_0, \ol\bv_0)$ in \eqref{eq:T}, and 
\begin{equation*}
\sigma(p,q)=
\begin{cases}
(3/q-3/p)/2+1/2 & \text{for}\,\,\,2\leq p\leq 3,\\
3/(2q) &  \text{for}\,\,\, 3<p <\infty.
\end{cases}
\end{equation*}

\subsection{Bound of $(\theta^1,\bw^1)$}
As $(\theta^1,\bw^1)$ satisfies the linear problem \eqref{eq:T} with the initial data $(\theta_0,\bw_0),$ we infer from \cite[Theorem 2.7]{EvBS2014} and \eqref{es:I_0_2} that
 \begin{multline}\label{es:thetaw1_1}
\sum_{\ell=1,2} 
\Big( \sup_{0<t<\min\{2,T\}} \big( \|e^{t} \theta^1\|_{H^{1}_{q_{\ell}}(\Omega)} 
+\|e^{t} \bw^1\|_{B^{2-2/p}_{q_{\ell},p}(\Omega)} \big)
+\|e^t (\theta^1,\bw^1)\|_{L_p(0,\min\{2,T\};H^{1,2}_{q_\ell}(\Omega))} \Big)\\
\lesssim \sum_{\ell=1,2} 
\Big( \|\theta_0\|_{H^1_{q_\ell}(\Omega)} 
+ \|\bw_0\|_{B^{2-2/p}_{q_\ell,p}(\Omega)} \Big)
\lesssim  \CI_0 +\CE_T^2 +\CE_T^3.
\end{multline}

In what follows, we suppose that $T>2.$
Firstly, the bounds in \eqref{es:sg_2} yield that 
\begin{equation}\label{es:thetaw1_2}
\begin{aligned}
\|(\theta^1,\bw^1)(t)\|_{L_{q_1}(\Omega)} 
&\lesssim t^{-3/(2q_1)} \vertiii{(\theta_0,\bw_0)}_{q_1,q_1/2}, \\
\|\nabla (\theta^1,\bw^1)(t)\|_{L_{q_1}(\Omega)} 
&\lesssim t^{- 3/(2q_1)-1/2} \vertiii{(\theta_0,\bw_0)}_{q_1,q_1/2}, \\
\|\nabla^2 \bw^1(t)\|_{L_{q_1}(\Omega)} 
&\lesssim t^{- 3/q_1} \vertiii{(\theta_0,\bw_0)}_{q_1,q_1/2}, \\
\|(\theta^1,\bw^1)(t)\|_{L_{q_2}(\Omega)} 
&\lesssim t^{-3(1/q_1-1/(2q_2))} \vertiii{(\theta_0,\bw_0)}_{q_2,q_1/2}, \\
\|(\nabla \theta^1,\nabla \bw^1,\nabla^2 \bw^1)(t)\|_{L_{q_2}(\Omega)} 
&\lesssim t^{- 3/q_1} \vertiii{(\theta_0,\bw_0)}_{q_2,q_1/2},
\end{aligned}
\end{equation}
for $1<t< T.$ Note the fact that 
\begin{equation*}
\sum_{\ell=1,2} \vertiii{(\theta_0,\bw_0)}_{q_\ell,q_1/2}
 \lesssim \CI_0''\lesssim \CI_0 +\CE_T^2 +\CE_T^3.
\end{equation*}
Then we combine the bounds \eqref{es:thetaw1_1} and \eqref{es:thetaw1_2} to obtain that 
\begin{equation}\label{es:thetaw1_3}
\begin{aligned}
 &\|\Jt^{3/(2q_1)}(\theta^1, \bw^1)\|_{L_\infty(0, T; L_{q_1}(\Omega))}
  +\|\Jt^{b} \nabla \theta^1\|_{L_p(0,T;L_{q_1}(\Omega))}
 +\|\Jt^{b}\nabla \bw^1\|_{L_p(0,T;H^{1}_{q_1}(\Omega))}\\
&+ \|\Jt^b(\theta^1, \bw^1)\|_{L_\infty(0, T; L_{q_2}(\Omega))} 
+ \|\Jt^{b}(\theta^1,\bw^1)\|_{L_p(0,T; H^{1,2}_{q_2}(\Omega))}
\lesssim \CI_0 +\CE_T^2 +\CE_T^3
\end{aligned}
\end{equation}
for the indices $b,q_1,q_2$ fulfilling
\begin{equation*}
\begin{aligned}
2<q_1 < 3< q_2<\infty,\quad
\frac{1}{q_1}=\frac{1}{3}+\frac{1}{q_2},
\quad bp'>1,\quad
\big(\frac{3}{2q_1}+\frac{1}{2}-b\big)p 
=\big(\frac{3}{q_1}- \frac{3}{2q_2}-b\big)p> 1.
\end{aligned}
\end{equation*}
Furthermore, using the equations of $(\theta^1, \bw^1)$ and \eqref{es:thetaw1_3}, we have 
\begin{equation}\label{es:thetaw1_4}
\begin{aligned}
&\sum_{\ell=1,2} \|\Jt^{b} \pd_t(\theta^1, \bw^1)\|_{L_p(0,T;H^{1,0}_{q_\ell}(\Omega))}\\
\lesssim & \sum_{\ell=1,2} \Big( 
\|\Jt^{b} \nabla \theta^1\|_{L_p(0,T;L_{q_\ell}(\Omega))}
+ \|\Jt^{b} \nabla \bw^1\|_{L_p(0,T;H^{1}_{q_\ell}(\Omega))} \Big)\\
 \lesssim & \,\,  \CI_0 +\CE_T^2 +\CE_T^3.
\end{aligned}
\end{equation}
Thus \eqref{es:thetaw1_3} and \eqref{es:thetaw1_4} furnish that 
\begin{equation}\label{main:tw1}
\CE_T (\theta^1,\bw^1) \lesssim   \CI_0 +\CE_T^2 +\CE_T^3.
\end{equation}

\section{Estimates of $(\theta^2,\bw^2)$}
\label{sec:decay}
In this section, we derive the decay estimates of
\begin{equation*}
\begin{aligned}
  (\theta^2,\bw^2)(\cdot,t)= &\lambda_0 \int_0^t T(t-s)(\rho^2,\bv^2)(\cdot,s)\,ds \\
  \end{aligned}
\end{equation*}
for $T(t)$ and $(\rho^2,\bv^2)$ given in the last section.
That is, $(\theta^2,\bw^2)$ satisfies the linear problem
\begin{equation} \label{eq:thw2}
	\left\{ \begin{aligned}
&\pd_t\theta^2 + \gamma_1 \, \dv\bw^2=  \lambda_0 \rho^2 
&&\quad&\text{in} &\quad \Omega \times (0,T), \\
&\gamma_1 \pd_t\bw^2 -\DV\big(\bS(\bw^2) - \gamma_2 \theta^2 \bI\big)
=\lambda_0 \bv^2 
&&\quad&\text{in}& \quad \Omega  \times (0,T), \\
&\big(\bS(\bw^2) - \gamma_2 \theta^2 \bI\big)\bn_{\Gamma} =0
&&\quad&\text{on}& \quad \Gamma \times (0,T), \\
&(\theta^2, \bw^2)|_{t=0} =(0,0)
&&\quad&\text{in}& \quad \Omega.
\end{aligned}
\right.
\end{equation}
In view of \eqref{es:rhov_2} and \eqref{es:rhov_12}, we set 
\footnote{We still keep the convention \eqref{def:Hmn} in this section.}
\begin{equation*}
\begin{aligned}
Y_T(\rho^2,\bv^2) 
&=\|\Jt^{b}(\rho^2,\bv^2)\|_{L_{\infty}(0,T;H^{1,0}_{q_1}(\Omega))}+\sum_{q\in \{q_1/2,\, q_1, q_2\}}  
\CM^{b}_{p,q}(T;\rho^2,\bv^2)
\lesssim \CI_0 + \CE^2_T + \CE^3_T.
\end{aligned}
\end{equation*}
Now, the main result of this section is stated as follows. 
\begin{thm}\label{thm:thetaw}
Let $b>0,$ $1<p<\infty$ and $2<q_1<3 <q_2<\infty$ satisfying the conditions:
\begin{equation} \label{cdt:thetaw}
 \frac{1}{q_1}=\frac{1}{3}+\frac{1}{q_2},
\quad bp'>1,\quad
\frac{3}{2q_1}+\frac{1}{2} -\frac{1}{p}>b\geq \frac{3}{2q_1}\cdot
\end{equation}
If $Y_T=Y_T(\rho^2,\bv^2)$ is finite, then 
\eqref{eq:thw2} admits a unique solution 
$(\theta^2,\bw^2).$ 
Moreover, there exists a positive constant $C$ such that
\begin{multline}\label{decay:1}
\sup_{0<t<T} \Jt^{3/(2q_1)} \|(\theta^2,\bw^2)\|_{L_{q_1}(\Omega)} 
+\|\Jt^b \nabla  \theta^2\|_{L_p(0,T;L_{q_1}(\Omega))}
+\|\Jt^b \nabla  \bw^2\|_{L_p(0,T;H^{1}_{q_1}(\Omega))} \\ 
+\|\Jt^b (\theta^2,\bw^2)\|_{L_p(0,T;H^{1,2}_{q_2}(\Omega))} 
+\sum_{\ell=1,2} \Big(\|\Jt^{b}(\pd_t\theta^2,\pd_t\bw^2)\|_{L_p(0,T; H^{1,0}_{q_\ell}(\Omega))}  \Big)
\leq C Y_T.
\end{multline}
\end{thm}

In the following, we only verify the decay estimates of spacial derivatives of $(\theta,\bw),$ namely,
\begin{multline}\label{es:goal}
\sup_{0<t<T} \Jt^{3/(2q_1)} \|(\theta^2,\bw^2)\|_{L_{q_1}(\Omega)} 
+\|\Jt^b \nabla  \theta^2\|_{L_p(0,T;L_{q_1}(\Omega))} 
+\|\Jt^b \nabla \bw^2\|_{L_p(0,T;H^{1}_{q_1}(\Omega))}\\
+\|\Jt^b (\theta^2,\bw^2)\|_{L_p(0,T;H^{1,2}_{q_2}(\Omega))} \lesssim Y_T.
\end{multline}
By admitting decay estimate \eqref{es:goal} for a while, the equations in \eqref{eq:thw2} imply that
\begin{equation*}
\sum_{\ell=1,2} \Big(\|\Jt^{b}\pd_t\theta\|_{L_p(0,T; H^1_{q_\ell}(\Omega))} 
 +\|\Jt^{b}\pd_t\bw \|_{L_p(0,T; L_{q_\ell}(\Omega))} \Big)\lesssim Y_T.
\end{equation*}

On the other hand, the last inequality in \eqref{es:sg_1} yields that 
\begin{multline}\label{es:thetaw_1}
\sup_{0<t<\min\{2,T\}} e^{t} 
\|(\theta^2,\bw^2)\|_{H^{1,2}_{q_1}(\Omega)} 
+\sum_{\ell=1,2}\|\Jt^b (\theta^2,\bw^2)\|_{L_p(0,\min\{2,T\};H^{1,2}_{q_\ell}(\Omega))}\\
\lesssim \sum_{\ell=1,2}\|\Jt^b (\rho^2,\bv^2)\|_{L_p(0,T;H^{1,2}_{q_\ell}(\Omega))}
\lesssim Y_T.
\end{multline}
In other words, \eqref{es:thetaw_1} gives the boundedness of $(\theta^2,\bw^2)$ in short time interval $(0,\min\{2,T\}).$ 
Thus we assume that $T>2$ in what follows, and prove that 
\begin{equation}\label{claim:decay}
\sup_{2<t<T} \Jt^{3/(2q_1)} \|(\theta,\bw)\|_{L_{q_1}(\Omega)} 
+\|\Jt^b \nabla  (\theta,\bw)\|_{L_p(2,T;H^{0,1}_{q_1}(\Omega))} 
+\|\Jt^b (\theta,\bw)\|_{L_p(2,T;H^{1,2}_{q_2}(\Omega))} \lesssim Y_T.
\end{equation}
Then \eqref{es:thetaw_1} and \eqref{claim:decay} yield \eqref{es:goal}. This completes the proof of Theorem \ref{thm:thetaw}.
\medskip

In order to establish the bound \eqref{claim:decay},  we take advantage of the decomposition 
\begin{equation}\label{def:I_k}
  (\theta^2,\bw^2)=\lambda_0 \Big(\int_0^{t/2}
  +\int_{t/2}^{t-1}+\int_{t-1}^t\Big) 
  T(t-s)(\rho^2,\bv^2)(\cdot,s)\,ds
  =I_1(t) +I_2(t) +I_3(t),
\end{equation}
and then we study $I_1(t),$ $I_2(t)$ and $I_3(t)$ respectively. 
For simplicity, we denote 
\begin{equation}\label{def:J_k}
J_{k,0}(t)=\|I_k(t)\|_{L_{q_{1}}(\Omega)},\quad
J_{k,1}(t) =  \|\nabla I_k(t)\|_{H^{0,1}_{q_{1}}(\Omega)},\quad 
J_{k,2}(t) = \|I_k(t)\|_{H^{1,2}_{q_{2}}(\Omega)},
\end{equation}
with $k=1,2,3.$
Furthermore, 
we also make use of the fact that
\begin{equation*}
\sum_{\ell=1,2}\Big\| \Jt^{b}\vertiii{(\rho^2,\bv^2)(\cdot,t)}_{q_\ell,q_1/2} \Big\|_{L_p(0,T)} 
\leq Y_T
\end{equation*}
for $\vertiii{(\rho^2,\bv^2)}_{q_{\ell},q_1/2}
=\|(\rho^2,\bv^2)\|_{L_{q_1/2}(\Omega)}
+ \|(\rho^2,\bv^2)\|_{H^{1,2}_{q_{\ell}}(\Omega)}.$

\subsection{Bound of $I_1(t)$}
By our assumption $1/q_1=1/3+1/q_2,$ we set 
\begin{equation*}
\fa=\frac{3}{2q_1}+\frac{1}{2} 
=\frac{3}{2}\big(\frac{2}{q_1}-\frac{1}{q_2} \big) = 1+\frac{3}{2q_2}>1.
\end{equation*}
Then by using \eqref{es:sg_2} and H\"older inequalities, we have 
\begin{equation}\label{es:I1_1}
\begin{aligned}
J_{1,0}(t)& \lesssim 
\int_0^{t/2} (t-s)^{-\frac{3}{2q_1}}
\vertiii{(\rho^2,\bv^2)(\cdot,s)}_{q_1,q_1/2}\,ds\\
&\lesssim t^{-\frac{3}{2q_1}}
\Big\| \Js^{b}\vertiii{(\rho^2,\bv^2)(\cdot,s)}_{q_1,q_1/2} \Big\|_{L_p(0,t/2)} 
\Big( \int_0^{t/2} \Js^{-bp'} \,ds\Big)^{1/p'}\\
&\lesssim  t^{-\frac{3}{2q_1}} Y_T,
\end{aligned}
\end{equation}

\begin{equation}\label{es:I1_2}
\begin{aligned}
J_{1,\ell}(t)& \lesssim 
\int_0^{t/2} (t-s)^{-\fa}
\vertiii{(\rho^2,\bv^2)(\cdot,s)}_{q_\ell,q_1/2}\,ds\\
&\lesssim 
 t^{-\fa}
\Big\| \Js^{b}\vertiii{(\rho^2,\bv^2)(\cdot,s)}_{q_\ell,q_1/2} \Big\|_{L_p(0,\,t/2)} 
\Big( \int_0^{t/2} \Js^{-bp'} \,ds\Big)^{1/p'}\\
&\lesssim    t^{-\fa} Y_T,
\end{aligned}
\end{equation}
for $\ell =1,2,$ and $2\leq t\leq T.$ 
As $(\fa-b)p>1$ by \eqref{cdt:thetaw},  
we obtain from \eqref{es:I1_1} and \eqref{es:I1_2} that 
\begin{equation}\label{es:I1}
\begin{aligned}
\sup_{2<t<T} \Jt^{\frac{3}{2q_1}} \|I_1(t)\|_{L_{q_1}(\Omega)} 
+\|\Jt^b \nabla I_1(t)\|_{L_p(2,T;H^{0,1}_{q_1}(\Omega))} 
+\|\Jt^b I_1(t)\|_{L_p(2,T;H^{1,2}_{q_2}(\Omega))} 
 \lesssim Y_T.
\end{aligned}
\end{equation}

\medskip

\subsection{Bound of $I_2(t)$}
Now, let us consider the bound of $I_2(t).$
Firstly, we use the condition $bp'>1$ to get
\begin{equation}\label{es:I2_1}
\begin{aligned}
\Jt^{\frac{3}{2q_1}} J_{2,0}(t)& \lesssim 
\int_{t/2}^{t-1} \vertiii{(\rho^2,\bv^2)(\cdot,s)}_{q_1,q_1/2}\,ds\\
&\lesssim \Big\| \Js^{b}\vertiii{(\rho^2,\bv^2)(\cdot,s)}_{q_1,q_1/2} \Big\|_{L_p(t/2,t-1)} 
\Big( \int_0^{t} \Js^{-bp'} \,ds\Big)^{1/p'} \lesssim  Y_T,
\end{aligned}
\end{equation}

On the other hand, note that $\int_{1}^{\infty} t^{-\fa}\,dt<\infty$ as $\fa >1.$
Analogous to \eqref{es:I1_2}, we obtain from \eqref{es:sg_2} 
and  H\"older inequalities that
\begin{equation}\label{es:I2_2}
\begin{aligned}
\Jt^{b} J_{2,\ell}(t) & \lesssim 
\int_{t/2}^{t-1} (t-s)^{-\fa}
\Js^b\vertiii{(\rho^2,\bv^2)(\cdot,s)}_{q_\ell,q_1/2}\,ds\\
&\lesssim \Big( \int_{t/2}^{t-1} 
(t-s)^{-\fa} \,ds \Big)^{1/p'}
\Big( \int_{t/2}^{t-1} 
(t-s)^{-\fa}
\Js^{bp}\vertiii{(\rho^2,\bv^2)(\cdot,s)}_{q_\ell,q_1/2}^p\,ds \Big)^{1/p}\\
&\lesssim  \Big( \int_{t/2}^{t-1} 
(t-s)^{-\fa}
\Js^{bp}\vertiii{(\rho^2,\bv^2)(\cdot,s)}_{q_\ell,q_1/2}^p\,ds \Big)^{1/p}.
\end{aligned}
\end{equation}

Hence Fubini's theorem and \eqref{es:I2_1} yield that
\begin{equation}\label{es:I2}
\sup_{2<t<T} \Jt^{\frac{3}{2q_1}} \|I_2(t)\|_{L_{q_1}(\Omega)} 
+ \|\Jt^b \nabla I_2(t)\|_{L_p(2,T;H^{0,1}_{q_1}(\Omega))} 
+\|\Jt^b I_2(t)\|_{L_p(2,T;H^{1,2}_{q_2}(\Omega))}  \lesssim Y_T.
\end{equation}
\medskip

\subsection{Bound of $I_3(t)$}
Since $(\rho^2, \bv^2)(t, \cdot) \in D(A)$ for any $t \in (0, T)$ as follows
from \eqref{sft**}, 
to estimate the norm of $I_3(t)$,  we use  
\eqref{es:sg_1} and H\"older inequalities and obtain that
\allowdisplaybreaks
\begin{align*}
\Jt^{b} J_{3,0}(t) &\lesssim 
\int_{t-1}^t \Js^{b}\|(\rho^2,\bv^2)(\cdot,s)\|_{H^{1,2}_{q_1}(\Omega)} \,ds\\
&\lesssim \int_{t-1}^t \,ds\,
\|\Js^{b}(\rho^2,\bv^2)(\cdot,s)\|_{L_{\infty}(0,T;H^{1,2}_{q_1}(\Omega))} \lesssim Y_T ,\\
\Jt^{b} J_{3,\ell}(t)& \lesssim 
\int_{t-1}^t \Js^{b}\|(\rho^2,\bv^2)(\cdot,s)\|_{H^{1,2}_{q_\ell}(\Omega)} \,ds
\\
& \lesssim  \Big( \int_{t-1}^t \,ds \Big)^{1/p'}
\Big(\int_{t-1}^t \Js^{bp}
\|(\rho^2,\bv^2)(\cdot,s)\|_{H^{1,2}_{q_\ell}(\Omega)}^p \,ds\Big)^{1/p}\\
& \lesssim  \Big(\int_{t-1}^t \Js^{bp}
\|(\rho^2,\bv^2)(\cdot,s)\|_{H^{1,2}_{q_\ell}(\Omega)}^p \,ds\Big)^{1/p},
\end{align*}

which together with Fubini's theorem and $b\geq 3/(2q_1)$ gives us 
\begin{equation}\label{es:I_3_1}
\begin{aligned}
&\sup_{2<t<T} \Jt^{\frac{3}{2q_1}} \|I_3(t)\|_{L_{q_1}(\Omega)} 
+\|\Jt^b \nabla I_3(t)\|_{L_p(2,T;H^{0,1}_{q_1}(\Omega))} 
+\|\Jt^b I_3(t)\|_{L_p(2,T;H^{1,2}_{q_2}(\Omega))} \\
\lesssim & Y_T+\sum_{\ell=1,2}\Big(  \int_1^T\Js^{bp}
\|(\rho^2,\bv^2)(\cdot,s)\|_{H^{1,2}_{q_\ell}(\Omega)}^p  \int_{s}^{s+1} 
\,dt\,ds\Big)^{1/p} \lesssim Y_T.
\end{aligned}
\end{equation}
At last, combing \eqref{es:I1}, \eqref{es:I2} and \eqref{es:I_3_1} yields \eqref{es:goal}.

\subsection{Complement of the proof of Theorem \ref{thm:main_1}}
In this subsection, we discuss the local wellposedness issue of \eqref{eq:LL_CNS_1} (or \eqref{eq:LL_CNS_2} equivalently). 
By the argument in \cite{EvBS2014}, it is not hard to see the following theorem. 
\begin{thm}\label{thm:lw_1}
Let $\Omega$ be a $C^2$ exterior domain in $\BR^3$ and let $2<p<\infty$ and $3<q_2<\infty.$ 
Assume that $(\rho_0,\bv_0) \in E_{p,q_2}(\Omega)$  satisfying the compatibility condition 
\begin{equation*}
\Big(\bS(\bv_0) - \big(P(\rho_e+\rho_0)-P(\rho_e)\big)\bI\Big)\bn_{\Gamma} = 0,
\end{equation*}
and $\|\rho_0\|_{H^1_{q_2}(\Omega)} 
+ \|\bv_0\|_{B^{2-2/p}_{q_2,p}(\Omega)} \leq D_0$ for some 
$D_0>0.$ 
Then there exists a (small) $T>0$ such that
problem \eqref{eq:LL_CNS_1} admits a unique solution $(\eta, \bu)$ with
$\|\eta\|_{L_{\infty}(0,T;L_{\infty}(\Omega))} \leq \rho_e/4$ and
\begin{equation*}
\begin{aligned}
\eta \in H^1_p\big(0, T;H^1_{q_2}(\Omega)\big),\,\,\, 
\bu \in H^1_p\big(0, T;L_{q_2}(\Omega)^3\big)
 \cap L_p\big(0, T; H^2_{q_2}(\Omega)^3\big).
\end{aligned}
\end{equation*}
Moreover, there exists a constant $C>0$ so that  
\begin{equation*}
\|\eta\|_{H^1_p(0, T; H^1_{q_2}(\Omega))} 
+ \|\bu\|_{H^1_p(0, T; L_{q_2}(\Omega))} 
+ \|\bu\|_{L_p(0, T; H^2_{q_2}(\Omega))} \leq CD_0.
\end{equation*}
\end{thm}

Now let us prove Theorem \ref{thm:main_1}. Assume that $\CI_0 \leq \ep $ for some small $\ep>0.$ Then according to Theorem \ref{thm:lw_1}, the solution $(\eta,\bu)$ satisfies that
$$\|\eta\|_{L_\infty(0,T;L_{\infty}(\Omega))}  \leq \rho_e/2,\quad
\int^T_0\|\nabla\bu(\cdot, s)\|_{L_{\infty}(\Omega)}\,ds\leq 1/2.$$
As $(\eta,\bu)=(\rho^1,\bv^1)+(\rho^2,\bv^2)+(\theta^1,\bw^1)+(\theta^2,\bw^2),$ 
the bounds \eqref{es:rhov_12}, \eqref{main:tw1} and \eqref{decay:1} yield that 
\begin{equation*}
\CE_T \leq C(\ep+ \CE_T^2 +\CE_T^3) 
\end{equation*}
for $\CE_T=\CE_{T}(\eta,\bu).$ Then the continuity of $\CE_T$ and \cite[Lemma 5.2]{SSZ2020}  give that $\CE_T  \leq M\ep $ for some constant $M$ so long as $\ep$ small enough.
So the solution $(\eta,\bu)$ can be extended to any larger time interval 
 by the standard bootstrap argument.


\section*{Acknowledgement}
Y.S. is partially supported by Top Global University Project and 
JSPS Grant-in-aid for Scientific Research (A) 17H0109; X.Z. is partially supported by NSF of China under Grant 12101457 and the Fundamental Research Funds for the Central Universities.


\end{document}